\documentclass[12pt,a4paper]{amsart}
\usepackage{amssymb, amstext, amscd, amsmath, color}
\usepackage{graphicx}

\usepackage{url}

\usepackage{hhline}

\begin{document}

\newtheorem{theorem}{Theorem}[section]

\newtheorem{corollary}[theorem]{Corollary}

\newtheorem{proposition}[theorem]{Proposition}

\newtheorem{lemma}[theorem]{Lemma}

%
\theoremstyle{definition}
\newtheorem{rem}[theorem]{Remark}



\newtheorem{definition}[theorem]{Definition}
\newtheorem{note}[theorem]{Note}
\newtheorem{eg}[theorem]{Example}
\newcommand{\Prf}{\noindent\textbf{Proof.\ }}
\newcommand{\bx}{\hfill$\blacksquare$\medbreak}
\newcommand{\upbx}{\vspace{-2.5\baselineskip}\newline\hbox{}%
\hfill$\blacksquare$\newline\medbreak}
\newcommand{\eqbx}[1]{\medbreak\hfill\(\displaystyle #1\)\bx}

\newcommand{\FFock}{\mathcal{F}}
\newcommand{\kil}{\mathsf{k}}
\newcommand{\Hil}{\mathsf{H}}
\newcommand{\hil}{\mathsf{h}}
\newcommand{\Kil}{\mathsf{K}}
\newcommand{\Real}{\mathbb{R}}
\newcommand{\Rplus}{\Real_+}

%

\newcommand{\bC}{{\mathbb{C}}}
\newcommand{\bD}{{\mathbb{D}}}
\newcommand{\bN}{{\mathbb{N}}}
\newcommand{\bQ}{{\mathbb{Q}}}
\newcommand{\bR}{{\mathbb{R}}}
\newcommand{\bT}{{\mathbb{T}}}
\newcommand{\bX}{{\mathbb{X}}}
\newcommand{\bZ}{{\mathbb{Z}}}
\newcommand{\bH}{{\mathbb{H}}}
\newcommand{\BH}{{\B(\H)}}
\newcommand{\bsl}{\setminus}
\newcommand{\ca}{\mathrm{C}^*}
\newcommand{\cstar}{\mathrm{C}^*}
\newcommand{\cenv}{\mathrm{C}^*_{\text{env}}}
\newcommand{\rip}{\rangle}
\newcommand{\ol}{\overline}
\newcommand{\td}{\widetilde}
\newcommand{\wh}{\widehat}
\newcommand{\sot}{\textsc{sot}}
\newcommand{\wot}{\textsc{wot}}
\newcommand{\wotclos}[1]{\ol{#1}^{\textsc{wot}}}
 \newcommand{\A}{{\mathcal{A}}}
 \newcommand{\B}{{\mathcal{B}}}
 \newcommand{\C}{{\mathcal{C}}}
 \newcommand{\D}{{\mathcal{D}}}
 \newcommand{\E}{{\mathcal{E}}}
 \newcommand{\F}{{\mathcal{F}}}
 \newcommand{\G}{{\mathcal{G}}}
\renewcommand{\H}{{\mathcal{H}}}
 \newcommand{\I}{{\mathcal{I}}}
 \newcommand{\J}{{\mathcal{J}}}
 \newcommand{\K}{{\mathcal{K}}}
\renewcommand{\L}{{\mathcal{L}}}
 \newcommand{\M}{{\mathcal{M}}}
 \newcommand{\N}{{\mathcal{N}}}
\renewcommand{\O}{{\mathcal{O}}}
\renewcommand{\P}{{\mathcal{P}}}
 \newcommand{\Q}{{\mathcal{Q}}}
 \newcommand{\R}{{\mathcal{R}}}
\renewcommand{\S}{{\mathcal{S}}}
 \newcommand{\T}{{\mathcal{T}}}
 \newcommand{\U}{{\mathcal{U}}}
 \newcommand{\V}{{\mathcal{V}}}
 \newcommand{\W}{{\mathcal{W}}}
 \newcommand{\X}{{\mathcal{X}}}
 \newcommand{\Y}{{\mathcal{Y}}}
 \newcommand{\Z}{{\mathcal{Z}}}

\newcommand{\sgn}{\operatorname{sgn}}
\newcommand{\rank}{\operatorname{rank}}

 \title[Polynomials for crystal frameworks]{
Polynomials for crystal frameworks and the rigid unit mode spectrum}

\author[S. C. Power]{S. C. Power}
\address{Dept.\ Math.\ Stats.\\ Lancaster University\\
Lancaster LA1 4YF \\U.K. }
\email{s.power@lancaster.ac.uk}


\thanks{2000 {\it  Mathematics Subject Classification.}
52C25, 74N05, 47N50 \\
Key words and phrases:
Crystal framework, matrix-valued function, rigidity operator, crystal polynomial, rigid unit mode.}

\date{}
\maketitle

\begin{abstract}
To each discrete translationally periodic bar-joint framework $\C$ in $\bR^d$
we associate a matrix-valued function $\Phi_\C(z)$ defined on the $d$-torus. The  rigid unit mode spectrum $\Omega(\C)$ of $\C$ is defined in terms of the multi-phases of phase-periodic infinitesimal flexes and is shown to correspond to the singular
points of the function $z \to \rank \Phi_\C(z)$
and also to the set of wave vectors of harmonic excitations which have vanishing energy in the  long wavelength limit.
To a crystal framework  in Maxwell counting equilibrium, which corresponds to $\Phi_\C(z)$ being square, the determinant of $\Phi_\C(z)$ gives rise to a unique multi-variable polynomial $p_\C(z_1,\dots ,z_d)$. For  ideal zeolites the algebraic variety of zeros of $p_\C(z)$ on the $d$-torus coincides with the RUM spectrum.
The matrix function is related to other aspects of idealised framework rigidity and flexibility and in particular leads to an explicit formula for the number of supercell-periodic floppy modes.
In the case of certain zeolite frameworks in dimensions $2$ and $3$  direct proofs are given to show the maximal floppy mode property (order $N$). In particular this is the case for the cubic symmetry sodalite framework and some other idealised zeolites.
\end{abstract}




\section{Introduction}

Let $\C$ be a mathematical crystal framework, by which we mean a connected structure in the Euclidean space $\bR^d$ consisting of a set $\C_e$ of framework edges, representing bars or bonds, with a corresponding set $\C_v$ of framework points (vertices), representing joints or atoms,  such that $\C_e$ is periodic with respect to a discrete translation group $\T$ of isometries of
$\bR^d$, with $\T$ of full rank. We consider mainly the case  $d=2, 3$ together with the locally finite assumption that $\C_e$ is generated by the translations of
a finite set  of edges.
Such a geometric bar-joint framework $\C$ can serve as a model for the essential geometry of the disposition of atoms and bonds in a material crystal $\M$. In this case  the vertices  have  atomic identifiers, such as   H, He, Li, B, ... , and  the chosen edges may correspond just to the strong bonds. The identification of strongly bonded molecular units, such as SiO$_4$ and TiO$_6$, imply a polyhedral net structure for $\C$ and in particular  aluminosilicate crystals  and zeolites provide in this way a fascinating diversity of tetrahedral nets in which every vertex is shared by two tetrahedra.

Material scientists are interested in the manifestation and explanation of various forms of low energy oscillation and excitation modes.
Of particular interest are the
rigid unit modes (RUMs) in crystalline materials, the low energy (long wavelength) modes of oscillation related to the relative motions of rigid units, such as the SiO$_4$
tetrahedral units in quartz.
The wave vectors of these  modes are observed in neutron
scattering experiments and have been shown to correlate closely with those for the modes observed in computer simulations with periodic networks of rigid units.
In both the experimental measurements and in the  simulations the background mathematical model is classical lattice dynamics  and the rigid unit mode wave vectors are observed where phonon dispersion curves display markedly low energy.
There is now a considerable body of literature tabulating the rigid unit mode wave vectors of various crystals and it has become evident that the primary determinant in a material $\M$  is the geometric structure of an associated abstract framework $\C$. This was outlined in the seminal paper of Giddy, Dove,  Pawley and  Heine \cite{gid-et-al}.
See also Swainson and Dove \cite{swa-dov},
Hammond et al \cite{ham-dov-zeo1997}, \cite{hamdov-zeo1998}
and  Dove et al \cite{dov-exotic}.
This experimental work shows that the wave vectors of RUMs often lie along lines and planes in reciprocal space. However, for many materials the wave vectors also lie on more mysterious curved surfaces. See also the recent computer assisted analysis of  Wegner \cite{weg}.

In what follows we develop a mathematical theory of rigid unit modes in idealized crystal frameworks.
As we shall demonstrate, this is essentially a linear first order theory and  one can side-step  lattice dynamical formulations that relate to higher energy phonons and their dispersion curves. In fact in Definition \ref{d:RumDefn2} we define the  RUM spectrum $\Omega(\C)$ of  an idealised
crystal framework $\C$, with given translation group, as the set of multi-phases for which there exists a nonzero phase-periodic infinitesimal flex. This form of the spectrum was first given in Owen and Power \cite{owe-pow-crystal}
as a byproduct of the analysis of square-summable infinitesimal flexes.
Mapping the $d$-torus to the unit cube in $\bR^d$ by taking logarithms
gives the usual wave vector  parametrisation space for RUMs used by crystallographers.
The spectrum  $\Omega(\C)$ leads naturally to a definition of the RUM dimension  $\dim_{\rm rum}\C$, which takes integer values from $0$ to $d$ and which gives a measure of the infinitesimal flexibility of $\C$. In the interesting special case of a crystal framework in Maxwell counting equilibrium (see Definition 2), such as, for example, a tetrahedral net framework derived from
an idealised zeolite, the spectrum $\Omega(\C)$ is determined as the zero set of a multi-variable polynomial $p_\C(z_1,\dots ,z_d)$ defined on the $d$-torus. This polynomial may vanish identically, which corresponds to the case $\dim_{\rm rum}(\C) = d$, and for $d=2,3$ this is also known as "order N". (See Theorem \ref{t:LocalFloppyMode}.) This property occurs for example in the case of the cubic form sodalite framework $\C_{\rm SOD}$, as we prove below in Section 7 by infinitesimal analysis. Interestingly, Kapko et al \cite{kap-daw-riv-tre} have recently conducted a simulation analysis to determine the extent of this property in  idealized zeolites.

The infinitesimal flex perspective is useful for several reasons. Firstly it brings into play the fairly well-established theory of infinitesimal rigidity for finite bar-joint frameworks and this is of significance for local flexibility. On the other hand the consideration of general infinitesimal flexes in infinite bar-joint frameworks gives a route to understanding and predicting the appearance of  linear components (lines, planes, hyperplanes etc) observed experimentally in RUM wave vector sets. In addition,  the first order infinitesimal flexibility properties of a crystal framework $\C$ are implicit in the infinite rigidity matrix $R(\C)$ of $\C$, and for  phase-periodic flexibility
this simplifies to the consideration of a finite function matrix $\Phi_\C(z)$  defined on the $d$-torus. This matrix function, which we also refer to as the symbol function of $\C$ (borrowing terminology from
Hilbert space  operator theory) also arises naturally from square-summable flex perspectives (\cite{owe-pow-crystal}) and may be
a useful tool more generally. When the matrix is square the \textit{crystal polynomial} $p_\C(z)$ of $\C$ is defined as a natural normalisation of its determinant.

In the development we give definitions, theorems, proofs and illustrative examples all of which lie within a mathematical theory of infinite bar-joint frameworks.
While the focus is on rigidity and flexibility properties related to the disposition of the bonds, rather than their strengths, the theory has the potential to relate to applied analysis and simulations.

In Sections 2 and 3 we give  examples of crystal frameworks and various spaces of infinitesimal flexes. In Sections 4, 5 and 6 we define the matrix function $\Phi_\C(z)$, the RUM spectrum $\Omega(\C)$ (Definition 8), the RUM dimension (Definition 9) and the crystal polynomial $p_\C(z)$. Also we give connections between phase-periodic infinitesimal flexes, so-called periodic floppy modes, and low energy harmonic excitations. (For discussions of wave vectors see Dove \cite{dov-book}, the account of phonon modes in Section 6 below, and the remarks following Definition 8 in Section 5 where we define the wave vector of a phase-periodic infinitesimal flex.)
In particular the matrix function $\Phi_\C(z)$
features in a counting formula for the periodic floppy modes in an $n$-fold supercell.

The final sections give determinations of $\Omega(\C)$ for a range of examples and some proofs. In particular we give the novel example of a two-dimensional zeolite
whose floppy modes are of order $N$.

\bigskip

\section{Crystal frameworks: terminology and examples.}
Let $G=(V,E)$ be a simple graph, finite or countable, with vertices $V=\{v_1,v_2,\dots \},$ and  $E \subseteq V \times V$
a countable set of edges, and let $p_1, p_2, \dots $ be a sequence of points in the Euclidean space $\bR^d$, with $p_i\neq p_j$ if $(v_i,v_j)$ is an edge. Then the pair $(G,p)$, with $p=(p_1,p_2,\dots )$ is  said to be a \textit{bar-joint framework} in $\bR^d$ with framework points, or joints, $p_i$
and framework edges, or bars, given by the line segments $[p_i, p_j]$
between $p_i$ and $p_j$ when $(v_i, v_j)$ is an edge in $E$. In all our examples in fact, the
framework points are distinct.

An \textit{isometry} of $\bR^3$ is a distance-preserving map $T:\bR^3\to \bR^3$. A
\textit{full rank translation group} $\T$ is a set of translation isometries
$\{T_k:k\in \bZ^3\}$ with  $T_{k+l}=T_k+T_l
$ for all $k,l$, $T_k \neq {\rm Id}$ if $k\neq 0$, and
such that the three \textit{period vectors}
\[
a=T_{\gamma_1}0, \quad b=T_{\gamma_2}0,\quad  c = T_{\gamma_3}0,
\]
associated with the generators $\gamma_1=(1,0,0), \gamma_2=(0,1,0), \gamma_3=(0,0,1)$
of $\bZ^3$ are not coplanar. Full rank translation groups in $\bR^d$ are similarly defined.

The following definition of a crystal framework $\C$ follows the formalism of  Owen and Power \cite{owe-pow-crystal} and pairs a  bar-joint framework of crystallographic structure with a group $\T$ of its translation symmetries. The definition
brings into play  the periodic partitioning of the vertices  and edges of $\C$ by the $\T$-translates of a finite geometrical \textit{motif} of framework vertices and edges.

\begin{definition}\label{d:crystalframework}
A crystal framework $\C=(F_v, F_e, \T)$ in $\bR^d$, with full rank translation group $\T=\{T_k:k\in \bZ^d\}$
and motif $(F_v, F_e)$, is a countable bar-joint framework $(G,p)$ with framework vertices $ p_{\kappa,k}$, for $1\leq \kappa \leq t,  k\in \bZ^d$, such that

(i) $F_v$ is a finite set of framework vertices, $\{p_{\kappa,0}:1\leq \kappa \leq t\}$ in $\bR^d$, and $F_e$ is a finite set of framework edges,

(ii) for each $\kappa$ and $k$ the vertex $p_{\kappa,k}$ is the translate $T_kp_{\kappa,0}$,

(iii) the set $\C_v$ of framework vertices is the  union of the disjoint sets
$T_k(F_v)$ for $ k\in \bZ^d$,

(iv) the set $\C_e$ of framework edges is the  union of the disjoint sets
$T_k(F_e)$ for $k\in \bZ^d$.
\end{definition}

This definition contains all the ingredients necessary for
the definition of rigidity matrices and operators associated with the various forms of periodic infinitesimal flexes that we shall consider. 
We also refer to the framework vertices simply as the framework points.
Natural choices for $\T$ are maximal translation subgroups of the crystallographic (spatial) symmetry group, subgroups respecting preferred symmetry directions,
and subgroups corresponding to supercell periodicity. 

\begin{center}
\begin{figure}[h]
\centering
\includegraphics[width=6cm]{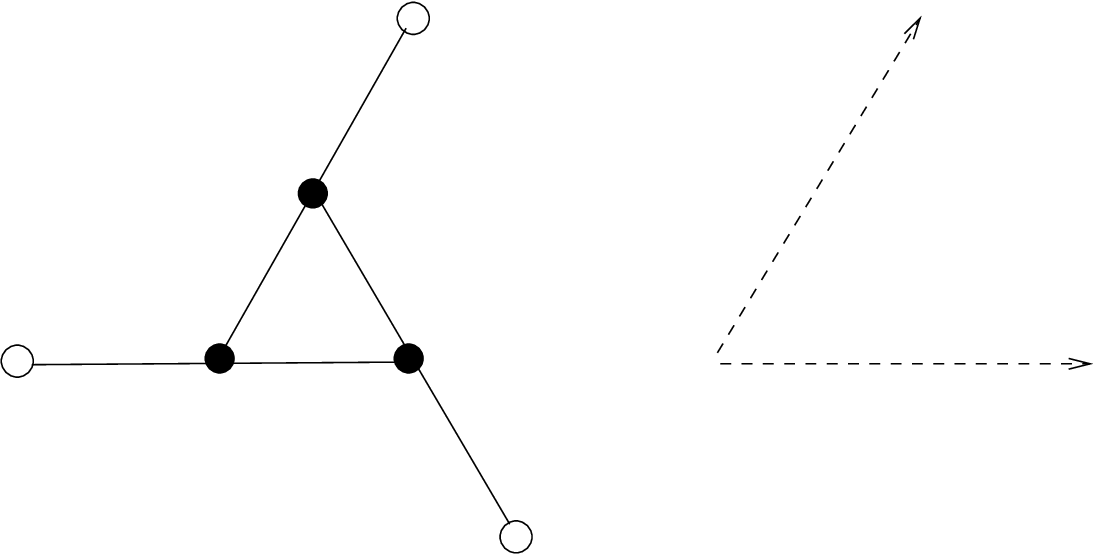}
\caption{Motif and period vectors for the kagome framework, $\C_{kag}$.}
\end{figure}
\end{center}

\begin{center}
\begin{figure}[h]
\centering
\includegraphics[width=3.5cm]{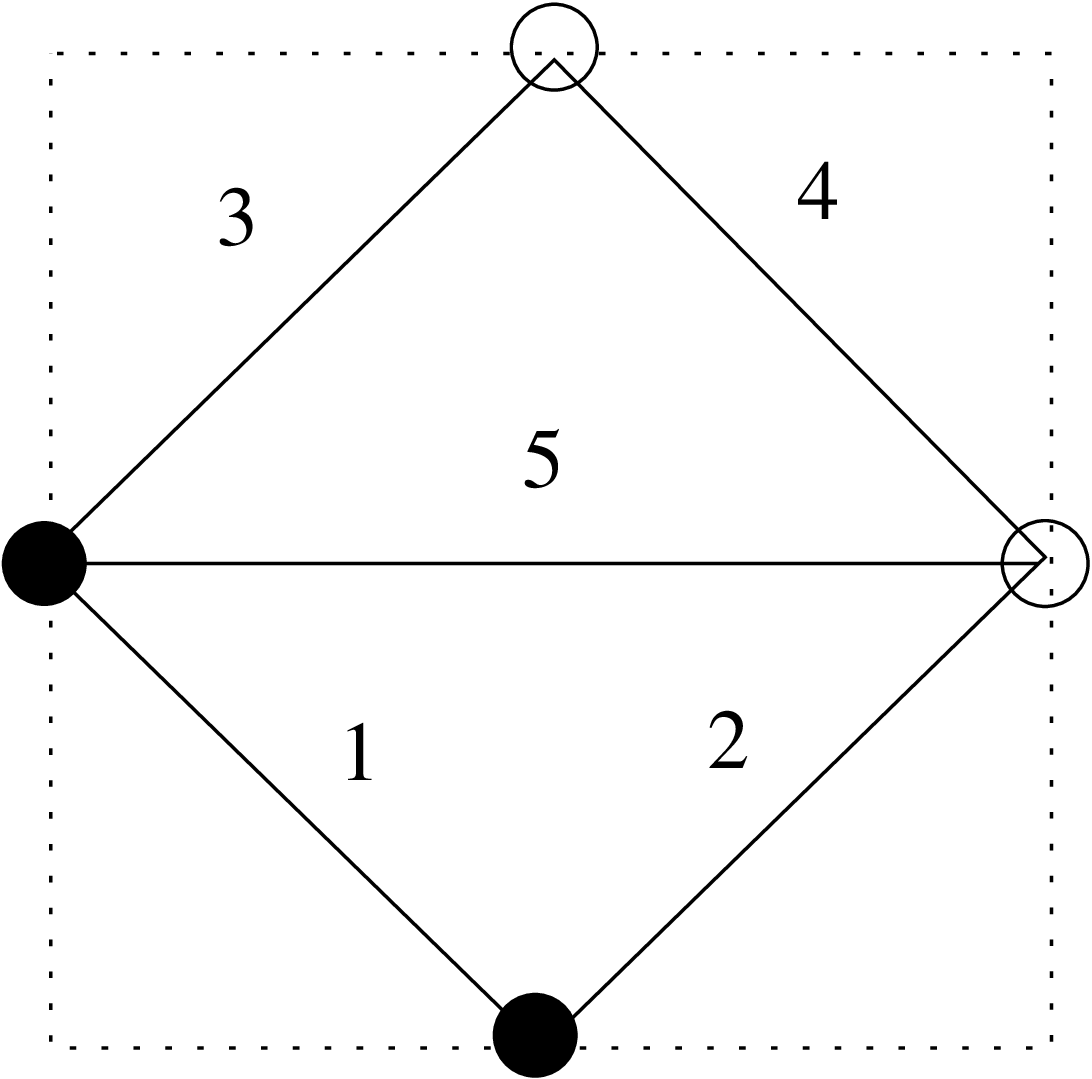}
\caption{A five-edged motif for the squares framework, $\C_{sq}$.}
\end{figure}
\end{center}


In Figures 1 and 2 motif choices are shown for the kagome framework $\C_{\rm kag}$ and the squares framework $\C_{\rm sq}$, where the filled vertices indicate the points of $F_v$ and
where the translation group is determined by the period vectors.
Thus $\C_{\rm kag}$ is the well-known  framework of pairwise corner-connected congruent equilateral triangles
in regular hexagonal arrangement, while $\C_{\rm sq}$ is a translationally periodic framework of corner-connected rigid square units.

One can similarly identify  motifs for other well-known frameworks and translation groups, such as
(i) the grid  framework $\C_{\bZ^d}$ in $\bR^d$ with $\C_v = \bZ^d$ and $\C_e$ equal to the set of line segments between nearest neighbours,  (ii) $\C_{\rm tri}$, the fully triangulated framework from the regular triangular tiling of the plane, and (iii) $\C_{\rm hex}$, the
crystal framework in the plane associated with the regular  hexagon
tiling.

In the examples we employ  a  mnemonic notational convenience with, typically,  $\C_{\rm xyz}$  with all lower case letters indicating a planar framework,  $\C_{\rm Xyx}$ indicating a 3D framework,
and $\C_{\rm XYZ}$ indicating a 3D framework which derives in a well-defined way from the zeolite with name XYZ.
For example we write $\C_{\rm Oct}$ to denote the basic regular octahedron net framework in three dimensions formed by corner-connected congruent octahedra with maximal cubic symmetry.
Also $\C_{\rm SOD}$, defined below,  derives from  the cubic form of the zeolite sodalite, while the companion framework $\C_{\rm RWY}$ derives from the sodalite RWY.
These conventions are useful, for example, when discussing subframeworks
lying in vector subspaces (slices).


The following definition is convenient.

\begin{definition}
A crystal framework $\C$ in $\bR^d$ is said to be in
\textit{Maxwell counting equilibrium}  if  $d|F_v|=|F_e|$
for some, and hence every, motif. If $d|F_v|<|F_e|$
then $\C$ is said to be \textit{edge rich} while if  $d|F_v|>|F_e|$
then $\C$ is said to be \textit{edge sparse}.
\end{definition}

We now define a number of illustrative crystal frameworks in dimensions 2 and 3 and in
Section 7 we compute their RUM spectra.
Of  particular interest with regard to rigidity and flexibility are the
$4$-regular ($4$-coordinated) frameworks in 2D and the $6$-regular frameworks in 3D, examples of which are provided by idealized zeolites in the sense of  Definition 3.

\medskip

\noindent{\bf Graphene and diamond bar-joint frameworks, $\C_{\rm gra}, \C_{\rm Dia},  \C_{\rm Dia}^2$.} The  usual visualisation of graphene is as a two-dimensional hexagonal bond-node network of carbon atoms with the geometry of $\C_{\rm hex}$. However $\C_{\rm hex}$ is edge sparse and this image is not suggestive of the strength of the material. If we view  the $C$-$C$-$C$ angles as rigid, or, equivalently, if we also view
\textit{second} nearest neighbours as  bonded, then
this leads to the edge rich crystal  framework, $\C_{\rm hex}^2$ say, implied by Figure 3. In the motif shown
we take two of the edges of one of the equilateral triangles to determine period vectors
and a corresponding translation group $\T$.

This crystal framework is of interest in its own right and  we also write it as $\C_{\rm gra}$ when viewed as a bar-joint framework in $\bR^2$. It may be assembled or decomposed in a number of ways to reveal substructure and in particular it may be constructed as a fusion of two congruent crystal subframeworks as follows.
Let $\C_{\rm tri}^+$ be obtained from the triangular framework $\C_{\rm tri}$ by adding bars, in which alternate triangles have three extra short bars added, connecting the triangle joints to a new joint at the centroid of the triangle.
Note that these added centroid joints are in natural one to one correspondence with the joints
of  $\C_{\rm tri}$ by a small translation.
Then $\C_{\rm gra}$ is congruent to the framework formed from the join of two copies of  $\C_{\rm tri}^+$, one of which is rotated
by $\pi$, and where the copies are connected by identifying the centroids of one copy with the non-centroid joints of the other, together with identification of the resulting
double edges.

\medskip

\begin{center}
\begin{figure}[h]
\centering
\includegraphics[width=7cm]{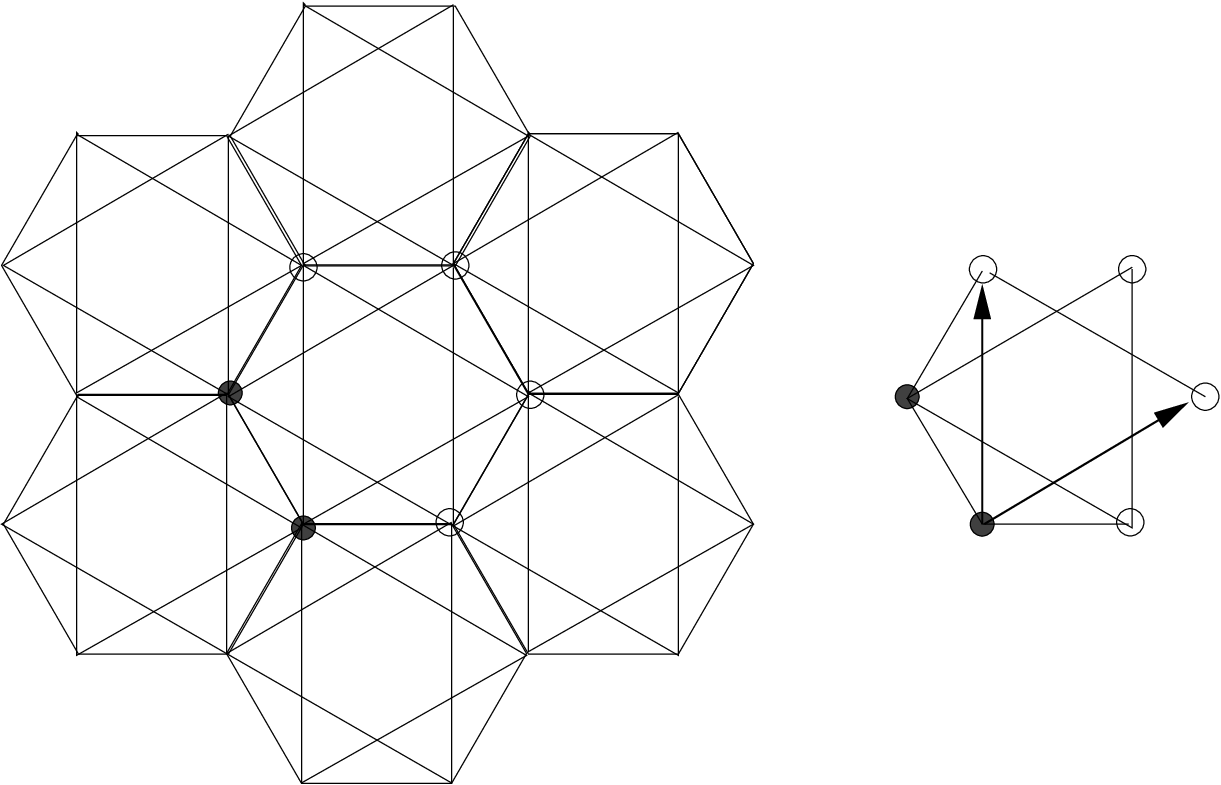}
\caption{Part of $\C_{\rm gra} = \C_{\rm hex}^2$ with choice of period vectors and motif.}
\end{figure}
\end{center}

\medskip

Similarly, crystalline diamond is usually indicated pictorially by a face-centred unit cell, with 14  C atoms at face centres and corners, plus 4 internal C atoms, and nearest neighbour connectivity.
Again, the implied 4-coordinated edge sparse bar-joint framework, $\C_{\rm Dia}$ say, does not of itself impart a sense of rigidity. It is natural for us to consider, once again, the derived first-and-second-nearest neighbour framework, and to take this as the definition of an associated bar-joint framework, which we denote $\C_{\rm Dia}^2$.
This too may be understood, or defined, in various constructive ways. For one such construction, echoing the graphene framework decomposition,  note that there is a bipartite red-blue colouring of the nodes of $\C_{\rm Dia}$ with face atoms red and internal  atoms blue say. The extra edges
of $\C_{\rm Dia}^2$ are either blue-blue or red-red. The red-red determined subframework we refer to as
the tetrahedron framework $\C_{\rm Tet}$.
Adding to this framework the blue-red edges of $\C_{\rm Dia}$ gives a framework we call $\C_{\rm Tet}^+$ (created by centroid addition).
It follows that $\C_{\rm Dia}^2$ is a join of two copies of $\C_{\rm Tet}^+$ (with reflected orientation), the join being effected by centroid/noncentroid identification, as before.

\bigskip

\noindent{\bf The cubic sodalite framework $\C_{{\rm SOD}}$.}
The crystal framework $\C_{{\rm SOD}}$ in three dimensions is built from 4-rings of tetrahedra
in a way which echoes the crystal structure of the cubic form of the zeolite sodalite. (See Figure 4.)

\medskip

The following general definition is convenient.
\begin{definition}
An ideal (or   mathematical)  zeolite in two (resp. three) dimensions is a
crystal framework $\C$ in the plane (resp. $\bR^3$) consisting of congruent triangles
(resp. congruent tetrahedra), each pair of which
intersect disjointly or at a common vertex and is such that every vertex is shared by two triangles (resp. tetrahedra).
\end{definition}

We remark that in databases \footnote{eg. http://www.iza-structure.org/databases/} material zeolite frameworks are most frequently indicated as a network of  "T atoms" corresponding to tetrahedral centres, each of which is $4$-coordinated with neighbouring $T$ atoms. This contrasts with the rigid unit view here of a tetrahedral net framework implied by the positions of  O atoms as vertices.

The 4-ring building units of  $\C_{{\rm SOD}}$ are oriented in the high symmetry arrangement indicated in Figure 4.
Six such rings may be placed on (the outside of the) six faces of an imaginary cube
so that the contact vertices sit on the midpoints of the edges of the cube.
This gives a finite bar-joint framework consisting of six regular 4-rings connected together to form a finite bar-joint framework which we call the sodalite cage framework.
With unit edge length for the tetrahedra the cube has
sidelength $1+\sqrt{2}$, while the three orthogonal period vectors (determining unit cell geometry) have length $2+\sqrt{2}$.

A motif for the framework can be given using the set $F_e$ of edges in three pairwise-connected
pairwise orthogonally oriented $4$-rings of the sodalite cage. The images of the edges of $F_e$ under the action of the associated isometry group $\T$
are essentially disjoint and generate the crystal framework $\C_{{\rm SOD}}$.
For an appropriate set $F_v$ an examination of the positioning of $F_v$ in the sodalite cage shows that one must take the vertices of $F_e$ except for $9$ redundant exterior vertices.

\begin{center}
\begin{figure}[h]
\centering
\includegraphics[width=6cm]{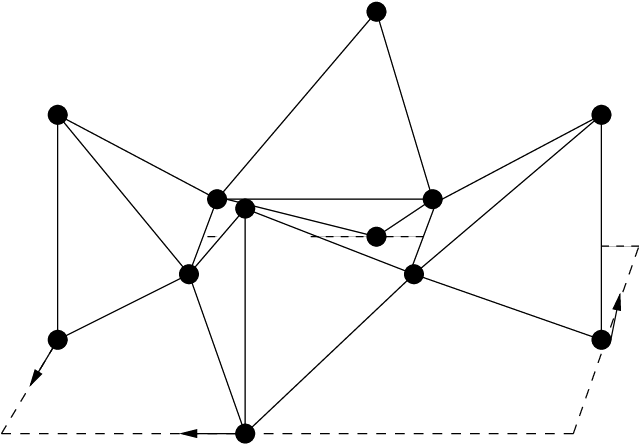}
\caption{The top 4-ring of the sodalite cage.}
\end{figure}
\end{center}

\bigskip

\noindent {\bf The kagome net framework $\C_{\rm Knet}$.}
We give two specifications of the \textit{kagome net} framework in three dimensions. Firstly, it may be constructed in a layered manner.  Form upward tetrahedral rigid unit frameworks on
alternate triangles of a two-dimensional kagome framework lying in the $xy$-plane. Similarly, form downward tetrahedra on the other triangles and thereby create a layer framework of pairwise connected tetrahedra. Parallel copies of such layers can be  joined at their exposed joints together to fill space,  creating, unambiguously,  the crystal framework we denote as $\C_{\rm Knet}$.

Alternatively, $\C_{\rm Knet}$ is a translationally periodic bar-joint framework with period vectors formed by three edges of a regular parallelapiped, with pairwise angles of $\pi/3$. Each parallelapiped contains two tetrahedral rigid units located at opposite "acute" corners of the parallelapiped  and with edge length half that of the parallelapiped edges. The  planar slices of $\C_{\rm Knet}$, determined by  each pair of period vectors, give   copies of $\C_{\rm kag}$.

\bigskip

\noindent {\bf The frameworks $\C_{\rm star}$ and $\C_{\rm oct}$.}
The kagome framework can be viewed as arising from the connection of translates of a regular 6-pointed star. There are analogous frameworks, $\C_{\rm star}$ and $\C_{\rm oct}$, arising from similar tilings using a regular 4-pointed star and an 8-pointed star respectively.
Figure 6, in the final section, indicates the (primitive) star template for $\C_{\rm star}$, while Figure 5 indicates tiling templates for four 2D zeolite crystal frameworks. More precisely, a natural choice of translation group for each of the associated bar-joint frameworks is that which is generated by horizontal and vertical translation. One can note that the third framework (with exterior angle $8\pi/12$), viewed as simply a bar-joint framework, is also recognisable as a congruent (rotated) copy of the bar-joint framework of  $\C_{\rm star}$. One can also confirm similarly that
the second and fourth frameworks are congruent by an isometry of $\bR^2$. The fourth framework here, with translation isometry group, is what we define as
the crystal framework   $\C_{\rm oct}$.

A motif $(F_v, F_e)$ for  $\C_{\rm oct}$ may be provided with $F_v$ the set of  four boundary vertices (indicated as solid vertices in the fourth template) plus the eight internal vertices (of the octagon), and
with $F_e$  the set of all 24 edges of the template.
Evidently, there is a smooth periodic edge-length-preserving continuous motion
(continuous or "finite" flex in the terminology of bar-joint frameworks) which "connects" these frameworks and which is parametrised by specification of the indicated exterior angle $\alpha$ say. This continuous motion, or evolution, maintains the squareness of the unit cells indicated in Figure 5 but evidently changes their edgelengths (and the period of translational periodicity). In Section 6 we give an indication of how the RUM spectrum evolves under this motion.
The motion here may be viewed as an example of the idealisation  of displacive phase transitions in materials (Dove \cite{dov-book}). We remark  that the derivative of this motion at any particular value of $\alpha$ gives a particular infinitesimal flex of the bar-joint framework associated with the  value $\alpha$. Such infinitesimal flexes are of \textit{affine type}
or \textit{flexible lattice type} and are not strictly periodic in the sense of Definition 6 below. For more on such infinitesimal flexes, which are associated with infinitesimal affine motions of the ambient space, see Borcea Streinu \cite{bor-str},
Power \cite{pow-aff} and Ross, Schulze and Whiteley \cite{ros-sch-whi}.

\begin{center}
\begin{figure}[h]
\centering
\includegraphics[width=9cm]{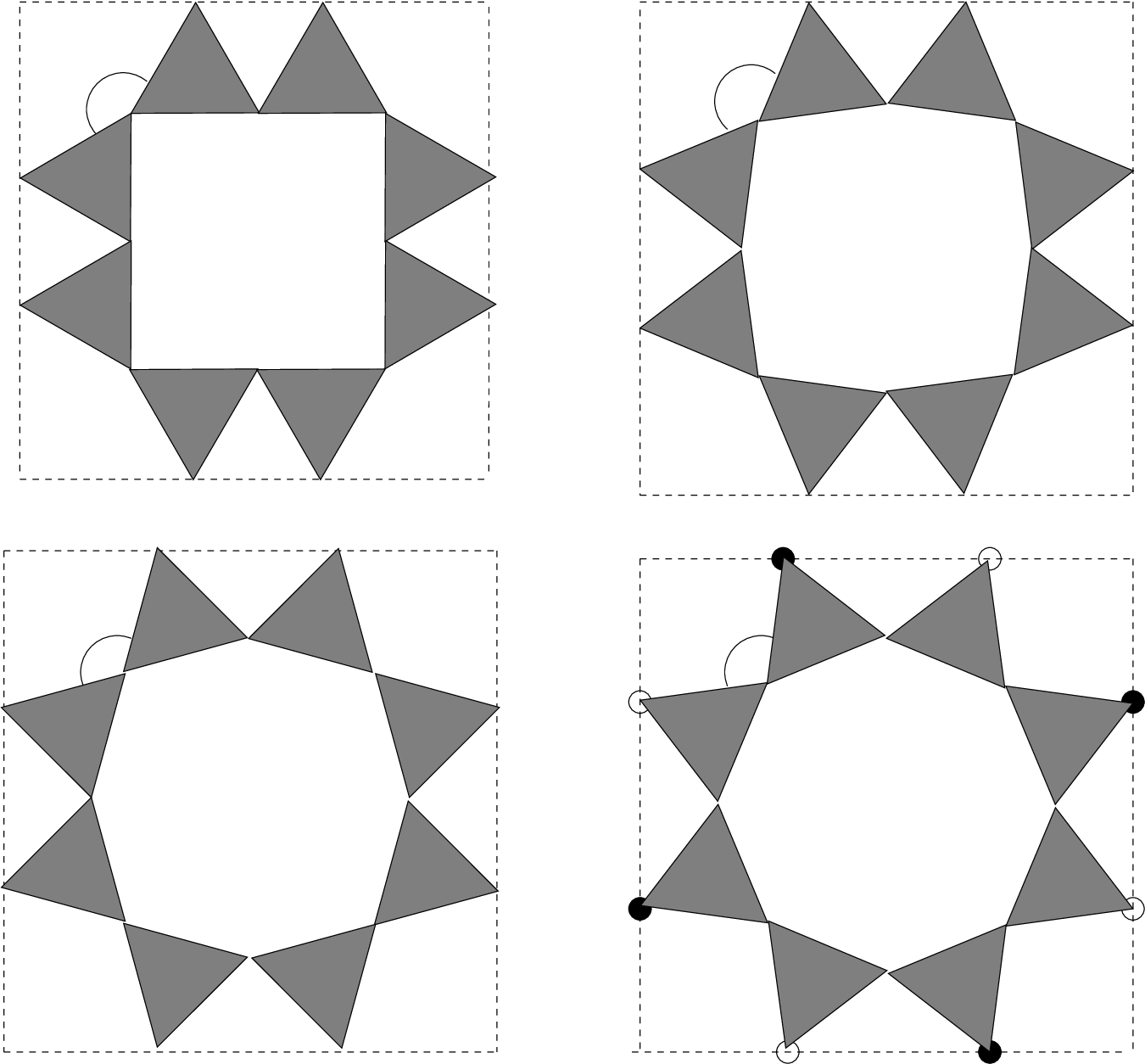}
\caption{Templates for 2D zeolite frameworks, with exterior angles
$10\pi/12, 9\pi/12, 8\pi/12, 7\pi/12$.}
\end{figure}
\end{center}

\noindent{\bf Further examples.}
Simple but informative examples of 3D zeolite frameworks can be built from  2D zeolite frameworks in various ways by layer constructions. With $\C_{\rm oct}$ for example, embedded in the $x, y$ plane of $\bR^3$, we may add bars and joints to obtain alternately upward and downward pointing tetrahedral units and so create a layer framework. These layers  may be joined consecutively at their exposed points to fill $\bR^3$ and thereby  create an associated ideal zeolite framework $\tilde{\C}_{\rm oct}$.
Similarly one can view $\C_{\rm Knet}$ as the framework $\tilde{\C}_{\rm kag}$.

We also note that interesting and diverse examples of mathematical crystal frameworks are implied by various tilings and periodic nets.
For an account of three-periodic nets and connections  with crystal chemistry see Delgado Friedrichs, O'Keeffe and Yaghi \cite{del-et-al-1}, \cite{del-et-al-2}. Such an (unlabelled) net, in any dimension, may be defined as a pair $(N,P)$, where $N$ the union of the edges of a crystal framework 
whose framework edges only intersect at framework vertices, and where $P$ is the set of framework points.

\bigskip

\section{Infinitesimal flexibility and rigidity.}
We now define various {flexes}  which act on the entire infinite crystal framework in a locally infinitesimal manner. The definition is the same as that for a finite bar-joint framework.

\begin{definition}\label{d:infflex}
An \textit{infinitesimal flex} of a finite or countable bar-joint framework $(G,p)$  is a vector $u=(u_i)$, with each component $u_i$ a vector  in $\bR^d$, such that for each edge $[p_i, p_j]$
\[
\langle p_i-p_j, u_i\rangle = \langle p_i-p_j, u_j\rangle.
\]
\end{definition}

Regarding the $u_i$ as velocity vectors this  asserts that for each edge the components of the endpoint velocities in the edge direction   are in agreement. This is equivalent to the assertion that  an infinitesimal flex is a velocity vector $v=(v_i)$ for which the distance deviation
\[
|p_i-p_j|-|(p_i+tu_i)-(p_j+tu_j)|
\]
of each edge is of order $t^2$ as the time parameter $t$ tends to zero.

We will not be concerned particularly with \textit{continuous flexes}, which are also called finite flexes (or finite edge-length-preserving deformations). For such flexes  each framework point undergoes a continuous motion $p_{\kappa,k}(t)$ such that edge lengths are preserved for all values of time $t$ in some range.
However we note that, as for  a finite framework, the derivative $u= p'(0)= (p_i'(0))$ of a continuous flex $p(t) = (p_i(t))$ with differentiable vertex trajectories
provides an infinitesimal flex $u$.

In the case of  a crystal framework in $\bR^d$
a {velocity vector}  is a doubly-indexed sequence $v$
of vectors $ v_{\kappa, k}$ in $\bR^d$ regarded as  instantaneous velocities applied to the frameworks vertices $p_{\kappa, k}$, and  it is convenient to consider
the vector space of all velocity sequences, written
as a direct product, namely
\[
\H_{atom}= \Pi_{\kappa,k}\bR^d.
\]
Thus, a real {infinitesimal flex} $u$ for the crystal framework $\C$ is a velocity vector $u$ in $\H_{atom}$ such that
\[
\langle   p_{\kappa, k} - p_{\tau, l}, u_{\kappa, k} - u_{\tau, l}\rangle  =0
\]
for each framework edge $[p_{\kappa, k}, p_{\tau, l}]$.
In particular the set of all infinitesimal flexes forms a vector subspace,  $\H_{\rm fl}$ say,
of $\H_{atom}$. Also each nontrivial infinitesimal isometry of $\bR^d$ gives rise to a one-dimensional vector subspace of $\H_{\rm fl}$.

The \textit{rigidity matrix} $R(\C)$ of $\C$ is a real infinite matrix defined as in the finite framework case.

\begin{definition}\label{d:RigidityMatrix}
The rigidity matrix $R(\C)$ of the crystal framework $\C$ in $\bR^3$ has rows
labelled by the edges $e= [p_{\kappa, k}, p_{\tau, l}]$ and
columns  labelled by the framework point coordinate indices
$(\kappa,x,k), (\kappa,y,k) , (\kappa,z,k)$. The  row for edge $e$ takes the form
\[
[\cdots 0 ~~ (p_{\kappa, k}-p_{\tau, l})~~ 0 \cdots 0 ~~(p_{\tau, l} -p_{\kappa, k})~~ 0 \cdots ]
\]
where the vector entry $(p_{\kappa, k}-p_{\tau, l})$ indicates that the three coordinates of this vector lie in the  columns for  $(\kappa,x,k), (\kappa,y,k) , (\kappa,z,k)$.
\end{definition}

The definition of $R(\C)$ for $d=2,4,5,\dots $, and also for general
countably infinite bar-joint frameworks (Owen and Power \cite{owe-pow-hon}) is essentially the same.
We remark that one may
take the view that $R(\C)$ is $1/2J(\C)$ where $J(\C)$ is the generalised Jacobian, evaluated at
the $p_{\kappa,k}$, for the infinite quadratic equation system
\[
|q_{\kappa,k}- q_{\tau,l}|^2 = d_e^2,
\]
where the  equations, labelled by the edges,  are in the coordinate variables of the points $q_{\kappa,k}$, and where the constants $d_e$ are the given lengths of the edges $e$ of $\C$.

It is natural to  consider various linear transformations that derive from $R(\C)$. To this end
let
$$\H_{bond} = \Pi_{e\in\C_e }\bR = \Pi_{e\in F_e,k\in\bZ^d} \bR
$$
be the space of real sequences $w=(w_{e,k})_{e\in F_e,k\in \bZ^d}$ labelled by the framework edges. Then $R(\C)$ gives a linear transformation $R:  \H_{atom}\to \H_{bond}$.
Indeed, each row of $R$ has at most $2d$ nonzero entries and the image $R(u)$ is given by the well-defined matrix multiplication $R(\C)u$. As for finite frameworks one has the following elementary proposition.

\begin{proposition}\label{p:infflex}
The infinitesimal flexes of the crystal framework $\C$ are the velocity vectors in $\H_{atom}$ that lie in the nullspace of the
linear transformation $R(\C)$.
\end{proposition}

Let us introduce notation for the natural basic sequences  of $\H_{atom}$ and $\H_{bond}$.
Write $\xi_x, \xi_y, \xi_z$ for the standard coordinate basis of $\bR^3$,
$\xi_x=(1,0,0)$ etc., and for $\sigma \in \{x,y,z\}$ write $\xi_{\kappa,\sigma, k}$ for the position indicator vector  in $\H_{atom}$ with
\[
(\xi_{\kappa,\sigma,k})_{\kappa',k'} = \delta_{\kappa,\kappa'}\delta_{k,k'}\xi_\sigma,
\]
where $\delta_{\kappa, \kappa'}$ is the Kronecker delta.
While $\H_{atom}$ does not have countable vector space dimension
its subspace of finitely nonzero sequences has the set $\{\xi_{\kappa,\sigma,k}\}$ as a vector space basis. However, the set is a generalised \textit{product type basis} for $\H_{atom}$ in the sense below.
In particular we may define the infinitesimal unit translation flex $u_x$ in the $x$ direction as the
well-defined infinite sum
\[
u_x=\sum_{\kappa,k} \xi_{\kappa,x,k}.
\]
Similarly we may write $\eta_{e,k}$  for the basic sequence in $\H_{bond}$ which is zero but for the value $1$ for the coordinate position $e,k$.

\medskip

Let $(G,p)$ be a countably infinite bar-joint framework in $\bR^d$.
A \textit{product type basis} for a subspace $\M$ of the velocity space $\H_{atom}$
of $(G,p)$ is a countable set $\S=\{w^1, w^2, \dots \}$ of  vectors in $\M$ such that,

(i) every vector $u$ in  $\M$ has a unique representation
\[
u = \sum_{n\in \bN} \alpha_n w^n,\quad \alpha_n\in \bR,
\]

(ii)  for each index $k$ only a finitely many elements $w^n$ of $\S$ have nonzero $k$th
component $w^n_{k}$.

The basic grid framework $\C_{\bZ^d}$  has evident nonzero infinitesimal flexes $u$ that act only on linear subframeworks  (copies of $\C_\bZ$ in $\bR^2$).
One can show that a set, $\S_d$  of representatives of all such flexes, is a product type basis for
$\H_{\rm fl}$.
In fact, we show elsewhere that it is possible to  identify product type bases for the vector space of all infinitesimal flexes for many other basic crystal frameworks.
Two examples are $\C_{\rm kag}$ and the 3D crystal framework $\C_{\rm Oct}$  for example.

\medskip

The following definition gives the context for the special classes of infinitesimal flexes of a crystal framework that  will concern us.

\begin{definition}
Let $\C$ be a crystal framework with translation group $\T$ as above.

(i) An infinitesimal flex (or  velocity sequence) $u$ is \textit{strictly periodic} if the following periodicity condition holds:
$u_{\kappa,k}=u_{\kappa,0}$ for all $k\in \bZ^d$.

\medskip
(ii) An infinitesimal flex (or  velocity sequence) $u$ is a \textit{supercell-periodic}  if $u_{\kappa,k}=u_{\kappa,0}$, for all $k$ in a subgroup $r_1\bZ \times \dots \times r_d\bZ$ for some positive integers $r_1,\dots , r_d$.

\medskip
(iii) An infinitesimal flex  $u$ is a \textit{local infinitesimal flex}
if $u_{\kappa,k}=0$ for all but finitely many values of $\kappa, k$.
\medskip

\end{definition}

Note the elementary fact that if there exists a local infinitesimal flex for $\C$ then this framework is rich in supercell-periodic flexes. Indeed if $u$ is such a local infinitesimal flex and  if $k\to \alpha_{k}$ is any supercell-periodic coefficient sequence
then the sum
\[
w = \sum_{k} \alpha_{k}T_ku,
\]
is a well-defined  supercell-periodic infinitesimal flex.

In Section 5 we turn attention to complex scalar infinitesimal flexes which are phase-periodic, the real and imaginary parts of which provide real infinitesimal flexes. It is such phase-periodic flexes that are closely allied to rigid unit mode wave vectors. They lead naturally to the formulation of a matrix-valued function associated with $\C$ and $\T$ and we
describe this association in the next section.

The strictly periodic infinitesimal flexes are also referred to as fixed lattice flexes. We remark that there is an interesting class of infinitesimal flexes  which lies outside our considerations here of (fixed lattice) rigid unit mode analysis, namely the affinely periodic infinitesimal flexes. Such "flexible lattice" flexes  allow, roughly speaking, an infinitesimal adjustment of the period vectors. (See also the comments on the frameworks in Figure 5.) For this and discussions of associated finite motions ee, for example, Borcea and Streinu \cite{bor-str}, Malestein and Theran \cite{mal-the},
Owen and Power \cite{owe-pow-crystal}, Power\cite{pow-aff}  and Ross et al \cite{ros-sch-whi}.
\medskip

\noindent {\bf Infinitesimal rigidity.} If a connected bar-joint framework $(G, p)$, finite or infinite, has no infinitesimal flexes other than rigid motion flexes then it is said to be \textit{infinitesimally rigid}. The simplest way in which this occurs  is when $(G, p)$ is \textit{sequentially infinitesimally rigid} (Owen and Power \cite{owe-pow-crystal}) in the sense that
it is the union of an increasing sequence of infinitesimally rigid finite frameworks. This is evidently the case for the edge rich frameworks $\C_{\rm tri}$, $\C_{\rm gra}$, $\C_{\rm Tet}$
and $\C_{\rm Dia}^2$. In particular it follows from the definitions below that the primitive RUM spectrum of each of these frameworks is trivial. Indeed, the RUM spectrum of a crystal framework is trivial when there are no phase-periodic infinitesimal flexes other than the strictly periodic flexes. For a sequentially rigid framework \textit{all} infinitesimal flexes are trivial rigid motion infinitesimal flexes, and the only phase-periodic flexes of this type are the strictly periodic rigid motion translation infinitesimal flexes.

On the other hand we remark that overconstrained
frameworks such as these edge-rich crystal frameworks are rich in periodic infinitesimal self-stresses. Following terminology for finite frameworks, a \textit{self-stress} $w = (w_e)_{e\in \C_e}$ of a crystal framework $\C$
is an assignment of scalars to edges such that for every framework point
$p_{\kappa, k}$ the finite vector sum
\[
\sum_{\tau,l:e=[p_{\kappa,k}, p_{\tau,l}]\in\C_e} w_e(p_{\kappa,k}-p_{\tau,l}),
\]
taken over all edges incident to $p_{\kappa, k}$, is equal to zero.
This is a companion notion to that of an infinitesimal flex and indeed $w$ is a self-stress if an only if $w$ lies in the nullspace of the
transpose of the rigidity matrix $R(\C)^T$. One may similarly consider subspaces of strictly periodic self-stresses and  phase-periodic self-stresses in a manner following the definitions for flexes.

One may relax the notion of infinitesimal rigidity to various forms of rigidity which are associated with a (possibly normed) space $S$ of velocity vectors. This is a viewpoint taken in Owen and Power \cite{owe-pow-hon}, \cite{owe-pow-crystal} leading to definitions of square-summable rigidity,
summable rigidity, vanishing flex rigidity ($c_0$-rigidity) and local flex rigidity ($c_{00}$-rigidity). 
It is an interesting issue to determine which classes of crystal frameworks are rigid in such a relative sense. In Kitson and Power
\cite{kit-pow} we analyse $c_0$-rigidity and its distinction from the other forms
of rigidity.

\section{The matrix function $\Phi_\C(z)$.}
A matrix-valued function, or \textit{symbol function}, for $\C$  is determined by the periodicity group $\T$ and the given motif $(F_v, F_e)$ as follows.

Write
$z=(z_1, \dots ,z_d)$, with $z_i\in \bC, |z_i|=1,$ to denote general points in the $d$-torus $\bT^d$.
Also, write $z^k$ for the monomial function
$z \to z^k$ from $\bT^d$ to $\bC$. Since ${z_i}^{-k}=\ol{z_i}^{k}$ for points on the circle $\bT$
we may think of general monomials $z^k$ as  products of the ${z_i}$ or $\ol{z}_i$ with just
non-negative powers.

It is convenient  to define the edge vector $v_e$ of the directed edge $e=[p_{\kappa,k}, p_{\tau,l}]$
as $v_e=p_{\kappa,k}- p_{\tau,l}$ and to write $v_{e,\sigma}$, for $1\leq \sigma \leq d$, for the coordinates of $v_e$.

\begin{definition}\label{d:Phi2}
Let $\C$ be a crystal framework in $\bR^d$ with  motif sets
\[
F_v=\{p_{\kappa,0}:1\leq\kappa\leq |F_v|\}, \quad F_e=\{e_i:1\leq i \leq |F_e|\}.
\]
Then $\Phi_\C(z)$ is the matrix-valued function on $\bT^d$ with rows labelled by the edges $e
=[p_{\kappa,k},p_{\tau,l}]\in F_e$ and with columns labelled by pairs $\kappa, \sigma$. As a matrix
of scalar function the entries are given by
\[
(\Phi_\C(z))_{e,(\kappa,\sigma)} = v_{e,\sigma}\bar{z}^k,
\]
\[
(\Phi_\C(z))_{e,(\tau,\sigma)}  = -v_{e,\sigma}\bar{z}^l,
\]
if $\kappa \neq \tau$,
while for a reflexive edge, with $\kappa = \tau$,
\[
(\Phi_\C(z))_{e,(\kappa,\sigma)}  = v_{e,\sigma}(\bar{z}^k-\bar{z}^l).
\]
The other entries are equal to the zero function.
\end{definition}

Different motifs for $\T$ give matrix functions that are equivalent in a natural way.
Indeed, replacement of a motif  edge (resp. vertex) by a $\T$-equivalent one  results
in the multiplication of the appropriate row (resp. columns) by a monomial. Thus in general two motif matrix functions $\Phi(z)$ and $\Psi(z)$ satisfy the equation
\[
\Psi(z)=D_1(z)A\Phi(z)BD_2(z),
\]
where $D_1(z), D_2(z)$ are diagonal matrix functions with monomial functions on the diagonal and where
$A, B$ are permutation matrices, associated with edge and vertex relabelling.

The next two examples and those we consider later  occur in two and three dimension and in this case we simply write $(z,w)$ and $(z,w,u)$ respectively for general points of $\bT^2$ and $\bT^3$.

\medskip
\noindent {\bf Example (a).} The motif for  $\C_{sq}$ implied by Figure $2$ has  $F_v$ equal to the ordered set $\{(1/2,0),(0,1/2)\}$ and $F_e=\{e_1,\dots ,e_5\}$. Here the period vectors, given by the sides of the parallelogram unit cell, are scaled with unit length.  It follows that the matrix function for $\C_{\rm sq}$ is
{\small
$$ \frac{1}{2}\left[ \begin {array}{cccc}
1 &-1&-1&1\\
-1 &-1&\bar{z}&\bar{z}\\
\bar{w} &\bar{w}&-1&-1\\
-\bar{w} &\bar{w}&\bar{z}&-\bar{z}\\
0&0&-2+2\bar{z}&0
\end {array} \right].
$$
}
If the final row of $\Phi_{\C_{\rm sq}}(z,w)$ is deleted then one has the matrix function for the realisation of the
square grid framework, $\C_{\bZ^2}$, when rotated by $\pi/4$. This framework is in Maxwell counting equilibrium and so the matrix function is square and we may compute
\[
\det \Phi_{\C_{\bZ^2}}(z,w)= 4\, \left( \bar{w}-\bar{z} \right)  \left(\bar{w}\bar{z}-1 \right).
\]
We consider further matrix function analysis for this example in Section 7, Example (f).

\medskip
\noindent {\bf Example (b).} With a choice of  labeling for the motif in Figure 1, with period vectors of length one, the matrix  function $\Phi_{\rm kag}(z,w)$ of the kagome framework
$\C_{\rm kag}$ takes the form
given by
{\small
\[
\Phi_{kag}({z},{w}) =
{\frac{1}{4}
\left[ \begin {array}{cccccc} -2&0&2&0&0&0\\
\noalign{\medskip}0&0&1&-\sqrt{3}&-1&\sqrt {3}\\
\noalign{\medskip}-1&-\,\sqrt {3}&0&0&1&\sqrt{3}\\
\noalign{\medskip}2&0&-2\,{z}&0&0&0\\
\noalign{\medskip}0&0&-1&\sqrt {3}&\ol{z}\ol{w}&-\sqrt{3}\ol{z}{w}\\
\noalign{\medskip}\ol{w}&\sqrt{3}\ol{w}&0&0&-1&-\sqrt {3}\end {array}
\right].}
\]
}
In this case the determinant is equal to a constant  multiple of
\[
\ol{z}\ol{w}{\left(z-1\right)  \left(w-1 \right)\left( z-w\right) }.
\]
For a different motif for the given translation group this determinant would change by a monomial factor.
\bigskip

\noindent{\bf Polynomials for crystal frameworks.}
Let $\C$ be a crystal framework in $\bR^d$ with a given isometry group $\T$.
If $\C$ is in Maxwell counting equilibrium then we may form the polynomial
~$\det(\Phi_\C(z))$ of the matrix function associated with a particular motif. This is a polynomial
in the coordinate functions $z_i$   and their complex conjugates $\ol{z_i}$, and is possibly identically zero. In the nonzero case we remove dependence on the motif and formally define the \textit{crystal polynomial} $p_\C(z_1,\dots ,z_d)$, associated with the pair $\C, \T$ and a
lexicographic monomial ordering,
as the product $\alpha z^\gamma\det(\Phi_\C(z))$ where the multi-power $\gamma$ and the scalar $\alpha$ are chosen so that
\medskip

(i) $p_\C(z)$ is a  linear combination of nonnegative power monomials,
\[
p_\C(z)= \sum_{\alpha \in \bZ^d_+} a_\alpha z^\alpha,
\]

(ii)  $p_\C(z)$  has minimum total degree, and

(iii)  $p_\C(z)$ has leading monomial with coefficient $1$.
\medskip

It is natural to order monomials  lexicographically, so that, for example, the monomial function
$z_1^2z_2$ has higher multi-degree than $z_1z_2^3$. In this way one defines the leading term of a multivariable polynomial. (See also the discussion in Cox, Little and O'Shea \cite{cox-lit-osh} for example.)

It follows that the crystal polynomial for the kagome framework and the (primitive case) translation group, as above, is
\[p_{kag}(z,w)={\left(z-1\right)  \left(w-1 \right)\left( z-w\right), }
\]
with lexicographic order $z>w$.
Also, for the grid framework $\C_{\rm \bZ^2}$ and the  non-axial translation group given above we see from the form of the determinant that
\[
p_{\rm \bZ^2}(z,w)= (z-w)(zw-1).
\]
For the grid frameworks $\C_{\rm \bZ^d}$ it is in fact more natural to take
the standard axial translation group $\T$ and a minimal motif which consists of a single vertex and $d$ edges, one for each axial direction. For this pair $\C, \T$ the crystal polynomial
is simply
\[
(z_1-1)(z_2-1)\dots (z_d-1).
\]

\section{Rigid Unit Modes and $\Phi_\C(z)$.}

We first show how $\Phi_\C(z)$ arises as a family of matrices parametrised by
points $z$ in the $d$-torus where the matrix for $z=\omega$ determines the possible existence of infinitesimal flexes
which  are periodic modulo the multi-phase $\ol{\omega}$.

Let $\K_{atom}, \K_{bond}$ be the complex scalar versions of the vector spaces  $\H_{atom}, \H_{bond}$.
Write $\K^\omega_a$ for the complex vector subspace space of complex velocity
vectors $v=(v_{\kappa, k})$ such that $v_{\kappa, k}=\omega^kv_{\kappa, 0}$ for
$\kappa\in F_v, k\in \bZ^d$. This is a finite-dimensional subspace of
$\K_{atom}$ of dimension $d|F_v|$.

Similarly let $\K^\omega_b\subset \K_{bond}$ be the subspace of the bond vector space
of complex sequences $w= (w_e)_{e\in \C_e}$ which are phase-periodic in this way for the phase $\omega$. Note that
the rigidity  matrix $R(\C)$ provides a linear transformation
$R^\omega:\K^\omega_a \to \K^\omega_b$
by restriction.
Indeed, with $d=3$, let $\gamma_i, 1 \leq i \leq 3$, denote the
usual generators for $\bZ^3$ and let $W_i$ and $U_i$
be the shift transformations on $\K_{atom}$ and $\K_{bond}$ respectively, with
\[
W_i: \xi_{\kappa,\sigma,k} \to \xi_{\kappa,\sigma,k+\gamma_i},
\]
\[
U_i: \eta_{e,k} \to \eta_{e, k+\gamma_i}.
\]
Then we have the commutation relations
\[
R(\C)W_i = U_iR(\C),\quad 1 \leq i \leq 3,
\]
and the identities $W_iu = \ol{\omega_i}u$, for $u\in \K^\omega_a$,
and  $U_iv = \ol{\omega_i}v$, for $v\in \K^\omega_b$. Thus for $u$ in  $\K^\omega_a$,
\[
U_i(R(\C)u)=R(\C)(W_iu)=R(\C)(\ol{\omega_i}u) = \ol{\omega_i}R(\C)u
\]
and so $R(\C)u \in \K^\omega_b$.

Let $\{\xi_{\kappa,\sigma}:
\kappa\in F_v, \sigma \in\{x,y,z\}\}$ be the natural basis for the column vector space $\bC^{3|F_v|}$. Write
$\xi_{\kappa,\sigma}^\omega$ for the
displacement vectors in $\K^\omega_a$ which "extend" the basis elements
$\xi_{\kappa,\sigma} $.
Formally, in terms of Kronecker delta symbol, we have
\[
(\xi^\omega_{\kappa,\sigma})_{\kappa',k} = \delta_{\kappa,\kappa'}\omega^k\xi_{\kappa,\sigma}.
\]
Similarly let $\eta_e, e\in F_e$, be the standard basis for $\bC^{|F_e|}$ and
write $\eta_e^\omega, e\in F_e$, for the natural  associated basis for $\K^\omega_b$, with
\[
(\eta_e^\omega)_{e',k}= \omega^k\delta_{e,e'}.
\]

\begin{theorem}\label{p:PhiRoleOne}
Let $\C$ be a crystal framework in $\bR^d$ with matrix function $\Phi_\C(z)$ and let $\omega\in \bT^d$.
Then the scalar matrix $\Phi_\C(\ol{\omega})$  is the representing matrix for the linear transformation
$R^\omega: \K^\omega_{a}\to \K_{b}^\omega$ with respect to the natural bases
$\{\xi_{\kappa,\sigma}^\omega\}$  and $\{\eta_e^\omega\}$.
\end{theorem}

\begin{proof} Let $\tilde{u}$ be a velocity vector in $\K_a^\omega$ determined by $u \in \bC^{d|F_v|}$ as above.
Let $e$ in $F_e$ be an edge of the form $[p_{\kappa,k},p_{\tau,l}]$ and
let $\langle \cdot , \cdot \rangle$ denote the  bilinear form on $\bC^d$. Then, from the definition of the rigidity matrix $R(\C)$, the $(e,0)^{th}$ entry of $R(\C)\tilde{u}$ in $\K_b^\omega$ can be written  as
\[
(R(\C)\tilde{u})_{e,0}= \langle v_e, \tilde{u}_{\kappa,k}  \rangle + \langle -v_e, \tilde{u}_{\tau,l}  \rangle
\]
\[
= \langle v_e, \omega^ku_\kappa\rangle  + \langle -v_e,  \omega^{l}u_\tau  \rangle
\]
\[
=\langle {\omega}^k v_e, u_\kappa\rangle + \langle -{\omega}^{l}v_e,u_\tau \rangle.
\]
This agrees with $(\Phi_\C(\ol{\omega})u)_e$, both in the case $\kappa \neq \tau$ and in the reflexive case $\kappa=\tau$ and the theorem follows.
\end{proof}

In particular the strictly periodic (one-cell-periodic) (real or complex) infinitesimal flexes are determined by the (real or complex) vectors in the
nullspace of the real matrix $\Phi(1,\dots,1)$. This \textit{periodic rigidity matrix} has rows carrying entries from the vectors $v_e, -v_e$ in the  case of nonreflexive edges of the motif (in the sense of Definition 7), with reflexive edges contributing zero rows.

The terminology of the following definition is justified in the next section.

\begin{definition}\label{d:RumDefn2}
The rigid unit mode spectrum (RUM spectrum) of the crystal framework $\C$ in $\bR^d$, with translation group  $\T$, is the set $\Omega(\C)$ of points
$\omega =(\omega_1, \dots ,\omega_d)$ in $\bT^d$
for which there is a nonzero vector $u$ in $\K_a^\omega$ which is an infinitesimal flex for $\C$.
\end{definition}

We also define the \textit{rigid unit modes} themselves as the nonzero infinitesimal flexes that give rise to points in the RUM spectrum. The \textit{mode multiplicity function} as the integer-valued function defined on $\Omega(\C)$  by
$\mu(\omega) = \dim \ker R^\omega$.

Note that from the  theorem we have 
\[
\Omega(\C)=\{\omega\in \bT^d:\ker \Phi(\ol{\omega})\neq\{0\}\}.
\]
In particular, from a commutative algebra perspective this set can be viewed as a real or complex algebraic variety.

\medskip

Evidently the RUM spectrum is a construct of the crystal framework $\C=(F_v, F_e,\T)$ and the ordering of coordinates matches the ordering of the generators 
of the translation group $\T=\{\T_k:k\in \bZ^d\}$. 
In the case of motif change, under a fixed  translation group $\T=\{\T_k:k\in \bZ^d\}$), one has two logically distinct crystal frameworks,
$\C=(F_v, F_e,\T)$ and $\C'=(F_v', F_e',\T)$  with
formally distinct symbol functions, $\Phi(z)$ and $\Psi(z)$ say. However, our earlier observation relating  these functions shows that in this case $\Omega(\C)$
and $\Omega(\C')$ are \textit{identical} subsets of the $d$-torus.
\medskip

\noindent{\bf Remarks.} In the next section we make precise the connection between rigid unit modes as we have defined them above and the low energy phonon modes that are of interest to material scientists, and are referred to as RUMs.
The convention in material science is to indicate  the set of  \textit{reduced wave vectors} that arise with these modes, rather than indicating a set of multi-phases, as we are doing here for their mathematical infinitesimal flex counterparts, but the conventions are simply related.  

If a material RUM has a wave vector  ${\bf k}=
({\bf k_1},{\bf k_2},{\bf k_3})$ then it has a
multi-phase $\omega= (\omega_1, \omega_2.\omega_3)$ in $\bT^3$ 
obtained by exponentiating, with $\omega_i=\exp(2\pi i{\bf k_i})$.
The  \textit{reduced wave vector} for the RUM is the reduction modulo $1$ in each coordinate and is the point   ${\bf k'}=
({\bf k_1'},{\bf k_2'},{\bf k_3'})$ in the unit cube $[0,1)^3$. It is obtained by taking the (principal) logarithms of each coordinate of the multi-phase.

For a simple crystal framework $\C$ in two or three dimensions (see Example (f) in Section 7 for example)  the set of RUM wave vectors
often consists of the intersection of $[0,1)^d$ with a union
of a finite number of points, lines and planes (hyperplanes for $d>3$) which are \emph{determined by equations over $\bQ$}. Also, in interesting cases the RUM wave vectors may fill all of $[0,1)^d$, with $\Omega(\C)=\bT^d$.
In these cases we say that $\Omega(\C)$ is a \textit{standard} RUM spectrum. Otherwise, borrowing terminology
from Dove et al \cite{dov-exotic}, we shall say that the RUM spectrum
is  \textit{exotic}. This includes the case of curves or curved surfaces in the unit cube. (The author is not aware  of examples of  crystal frameworks whose RUM spectrum has isolated irrational points or "exposed" irrational lines.)
\medskip

\noindent {\bf The dimension $\dim_{\rm rum}(G,p)$.}
Returning to the RUM spectrum recall that  $\Omega(\C)$ is a well-defined set in $\bT^d$ determined by the underlying bar-joint framework $(G,p)$ and a translation group $\T=\{\T_k:k\in \bZ^d\}$, and where  generator permutations for $\T$ correspond to a coordinate
permutation. We now define the RUM dimension of $(G,p)$, which takes an integer value between $0$ and $d$ inclusive.

We first define the \textit{primitive RUM spectrum} 
$\Omega_{\rm prim}(G,p)$ of the crystallographic bar-joint framework $(G,p)$ as the RUM spectrum for the  crystal framework $\C$ associated with $(G,p)$ and a \emph{maximal} translation group of isometric automorphisms of $(G,p)$.
(The terminology borrows  from the notion of a primitive unit cell
for a crystallographic set in $\bR^d$.)  To see that $\Omega_{\rm prim}(G,p)$ is well-defined, up to permutation of the coordinates, 
we first recall the classical fact of Bieberbach that a crystallographic group in any number of dimensions has a unique maximal normal free abelian subgroup.  In our setting this entails that two maximal translation subgroups $\T=\{T_k:k\in \bZ^d\}$ and $\T'=\{T_k':k\in \bZ^d\}$ of the  isometric (spatial) automorphism subgroup are congruent by an isometry $Z$ of $\bR^d$ which effects an automorphism
$Z$ of $(G,p)$. In this case we have  $T_k'=ZT_kZ^{-1}$ for all $k$. Moreover, in view of our earlier discussion we may assume that the motif $(F_v, F_e)$ for $\T$ is given and that the motif $(F_v', F_e')$ for 
$\T'$ is chosen as the image of  $(F_v, F_e)$ under $Z$, with a corresponding labelling. It follows that the respective symbol functions $\Phi(z)$ and $\Phi(z)'$ are simply related. Indeed, let 
$S$ be the linear isometry component of $Z$. Then, in the notation for the symbol functions the new motif edge vector $v_e'$ associated with framework edge $Ze\in F_e'=ZF_e$  is the  vector $Zv_e= Zp_{\kappa,k}-Zp_{\tau,l}$ which, being a difference, is equal to $Sv_e$. It follows from this that
\[
\Phi(z)' = \Phi(z)\tilde{X}
\]
where $\tilde{X}$ is an invertible  block diagonal (scalar) 
matrix $X\oplus \dots \oplus X$ (with $d|F_v|$ summands ). The well-definedness of the primitive RUM spectrum of $(G,p)$ now follows.

\begin{definition}\label{rumdimension}
Let $(G,p)$ be the crystallographic bar-joint framework, that is, a bar joint framework that underlies a crystal framework. Then the RUM dimension  $\dim_{\rm rum}(G,p)$ of
$(G,p)$ is  the real dimension of the real  algebraic variety
$\Omega_{\rm prim}(G,p)$.
\end{definition}

The dimension here can be considered as the  topological dimension of
the manifold of nonsingular points in case $\Omega_{\rm prim}(G,p)$ is irreducible. Otherwise it is the maximal such dimension over irreducible components. 
In fact we see below that the dimension of the RUM spectrum $\Omega(\C)$ of a crystal framework does not depend on the choice of translation group in view of a simple relationship between the RUM spectrum and the primitive RUM spectrum. Thus we may view the RUM dimension of $\C$ as this common dimension.

In view of the determinations in Section 7 and our comments below we shall see that
\[
 \dim_{\rm rum}(\C_{\rm sq}) = 0, \dim_{\rm rum}(\C_{\rm star}) = 1,
\dim_{\rm rum}(\C_{\rm kag}) = 1, \dim_{\rm rum}(\C_{\rm oct}) = 1,
\]
in two dimensions, and in higher dimensions we have
\[
\dim_{\rm rum}(\C_{\bZ^d}) = d-1,\dim_{\rm rum}(\C_{\rm Knet}) =2, \dim_{\rm rum}(\C_{\rm Oct})=1 , \dim_{\rm rum}(\C_{\rm SOD})=3.
\]

\medskip

 For a framework in Maxwell counting equilibrium the variety $\Omega(\C)$ is simply the zero set of $p_\C(z)$. For the kagome framework, for example, the  polynomial is $(z-1)(w-1)(z-w)$ and we obtain the set which is the union of the three curves on $\bT^2$ defined by $z=1$, $w=1$ and $z=w$. In terms of wave vectors this translates to the union of the three parametrised lines $(0, \alpha)$, $(\alpha,0)$ and $(\alpha,\alpha)$. Thus the RUM dimension is $1$.

When $\C$ is edge rich, with $|F_e|>d|F_v|$ then one may instead form the finite family of polynomials of the $d|F_v|\times d|F_v|$ submatrices of $\Phi_\C(z)$. Then the RUM spectrum will be a variety contained in the intersections of the zero sets of these polynomials on the $d$-torus.

We remark that the RUM spectrum will generally carry symmetries reflecting the point group symmetries of the crystal framework. Even so the point group may be trivial and the following rather theoretical inverse problem  may well have an affirmative answer.

\medskip

\noindent \textit{Problem.} Let $q(z,w)$ be a polynomial with real coefficients with $q(1,1)=0$.
Is there a crystal polynomial $p(z,w)$ whose zero set on the $2$-torus is the same as that for $q(z,w)$ ?
\bigskip

\noindent {\bf Floppy modes and their asymptotic order.}
In applications the term \textit{floppy mode} often refers to rigid unit flexibility and oscillation within a large  supercell and there is interest in the asymptotic order of the number of such modes as the supercell
dimensions tend to infinity. In particular, a so-called \textit{order $N$} crystal structure  (to use terminology employed by material scientists) is one for which the asymptotic order agrees with the order of the number of atoms in the supercell, which is of order $N= n^3$ in an $n\times n\times n$
supercell of a 3D crystal.
We now formalise this terminology in the direction of infinitesimal flexes  and indicate connections with the RUM spectrum.

\begin{definition}\label{d:FloppyMode}
Let $\C$ be a crystal framework in $\bR^d$ with translation group $\T=\{T_k:k\in \bZ^d\}$.

(i) An $n$-fold periodic floppy mode of  $\C$ is a nonzero real vector $u=(u_{\kappa,k})$ in the nullspace (kernel) of the rigidity matrix $R(\C)$
which is periodic for the subgroup 
$(n\bZ)^d$. That is,
$u_{\kappa,k}=u_{\kappa,0}$ for all $k \in (n\bZ)^d$.

(ii) $\nu_n$ is the dimension of the real linear space of real $n$-fold periodic floppy modes.

(iii) A crystal framework $\C$ in $\bR^3$ is of \textit{order $N^\alpha$} for floppy modes, where $\alpha = 0, \frac{1}{3}, \frac{2}{3}$ or $1$, if  $\nu_n \geq Cn^{3\alpha}
$ for all $n$ for some $C>0$, while there is no such constant for the power $\alpha +\frac{1}{3}$.
In particular $\C$ is said to be of order $N$ if $\nu_n \geq cn^{3}$ for some constant
$c>0$.
\end{definition}

Since the real and imaginary parts of a complex infinitesimal flex are real infinitesimal flexes
it follows that $\nu_n\leq \dim \ker R_n(\C)\leq 2\nu_n$ where $R_n(\C)$ is the rigidity matrix for $n$-fold
periodicity viewed as a complex vector space linear transformation. Thus, in considerations of asymptotic order we may more conveniently consider the complex scalar case. The matrix $R_n(\C)$ is the rigidity matrix for strict periodicity relative to the subgroup $\T'$ of $\T$ corresponding to the index subgroup $n\bZ\times \dots \times n\bZ$. Accordingly it is given as the periodic rigidity matrix associated with
a motif for the $n$-fold supercell. Such a motif  can be taken simply as the union of $n^d$ translates of the given motif. More conveniently, it is possible to explicitly block diagonalise $R_n(\C)$ (as a complex vector space transformation)  as a direct sum (even an orthogonal direct sum for natural inner product) of the matrices $\Phi(\omega)$ as $\omega$ ranges over the set of points,
$\bT_n$ say, with coordinates $\omega_j$ of the form $e^{2\pi k_j/n}$, where $0\leq k_j< n$ are integers.
This then gives the following counting formula for floppy modes:
\[
\dim \ker R_n(\C) =
\sum_{0\leq k_i<n,1\leq i\leq d} \dim\ker(\Phi_\C(\omega^k))
\]
where $\omega^k=(e^{2\pi ik_1/n}, \dots , e^{2\pi ik_d/n})$. This formula resolves a question posed by Simon Guest. An elementary  direct proof  follows from the fact that nonzero vectors $u,v$ from distinct nullspaces $\ker \Phi_\C(\omega^k)$ are linearly independent, on the one hand, and that, on the other hand, by the usual averaging arguments, any $n$-fold periodic flex
may be decomposed as a sum of pure frequency $n$-fold periodic infinitesimal flexes. By "pure frequency" we mean phase-periodic in each coordinate for some $n^{th}$ root of unity (depending on the coordinate). Details are given in the Appendix.

It is of interest then to consider the rational subset of the RUM spectrum corresponding to periodic floppy modes, namely
\[
\Omega_{\rm rat}(\C) := \bigcup_{n=1}^\infty (\Omega(\C)\cap \bT^d_{n})
\]
and to ask:
\medskip

To what extent does the asymptotic order of the periodic floppy modes determine the RUM dimension ?

\medskip

In the case of crystal frameworks with a primitive RUM spectrum which is \emph{standard} in the above sense, there is in fact a close connection. We make this clear below in the case of order $N$ (the maximal order).
In the exotic case  one should expect examples
where the rational points of the  RUM spectrum are not dense.
It would be of theoretical interest to identify, for example, a curved RUM spectrum only containing a finite number of rational points. Possibly the regular octagon ring framework $\C_{\rm oct}$ has this property.

For the proof of Theorem 2 we note the following lemma.

\begin{lemma}
(i) Let $\C$ be a $d$-dimensional crystal framework with motif set $(F_v, F_e)$ and RUM spectrum $\Omega(\C)\subseteq
\bT^d$. Then
\[
d-1+|\Omega(\C)\cap\bT^d_n|\leq \dim \ker R_n(\C)\leq d|F_v||\Omega(\C)\cap\bT^d_n|
\]
where $|F_v|$ is the number of vertices in the partition unit cell and where
$\bT^d_n$ is the "discrete torus" (in the $d$-torus $\bT^n$) determined by $n^{th}$ roots of unity.

(ii) If $\dim \ker R_n(\C) \geq cn^\alpha$ for some $c>0, \alpha>0$, then $\dim(\Omega(\C))\geq \alpha$.

\end{lemma}

\begin{proof}
(i) The counting formula implies the second inequality
since $\dim \ker (\Phi_\C(\omega)) \leq d|F_v|$ for all $\omega$.
Also, if $\omega^k \in \Omega(\C)\cap \bT^d_n$ then $\dim \ker \Phi_\C(\omega^k)\geq 1$,
while for wave vector ${\bf k}=(0,0,0)$ we have \\
$\dim\ker(\Phi_\C(1, \dots , 1))\geq d$, since there are certainly $d$ linearly independent
translation infinitesimal flexes. Thus the first inequality follows.

(ii) follows  from (i) since for any algebraic variety
$\Omega$, if the dimension is less than the integer $\alpha$ then the cardinality
of $\Omega\cap \bT^d_n$ is at most of order $n^{\alpha-1}$.
\end{proof}

It can be shown by direct linear algebra, as we now sketch, that if a crystal framework has order $N$ then there exists a local infinitesimal flex.

\begin{theorem}\label{t:LocalFloppyMode}
With the notation above the following properties are equivalent for a crystal framework $\C$ in $\bR^d$.

(i) $\C$ has a local infinitesimal flex.

(ii) $\C$ is of order $N$.

(iii)  $\dim_{\rm rum}(\C)=d$.

(iv)  $\Omega(\C) = \bT^d$.
\end{theorem}
\begin{proof} 
To see that (i) implies (ii) note that if $u$ is a nonzero local infinitesimal flex and $\omega$ is any multi-phase in $\bT^d$
then the sum
\[
v=\sum_{k\in\bZ^d} \omega^kT_ku
\]
is a phase-periodic infinitesimal flex. Also it is nonzero for almost every $\omega$. Thus (iv) holds and hence (iii) and (ii).

If (ii) holds then (iii), and hence (iv), follows from the lemma and the fact that $\Omega(\C)$  is a real algebraic variety in $\bT^d$.

Since (iv) implies (ii) it remains to show that (ii) implies (i) and this we do in the Appendix.
\end{proof}
\medskip

\noindent {\bf RUM spectrum verses primitive RUM spectrum.}
Note that if one doubles all the period vectors for $\C$ to obtain $\C'$ then it follows  that the new RUM spectrum  contains the range of the old spectrum under the argument doubling map, $\pi: (w_1, w_2, w_3) \to (w_1^2, w_2^2, w_3^2)$.
This follows immediately from the definition. The new symbol function, the number of  rows and columns of
which have increased $2^d$-fold, is less useful at this point.
In fact the map $\pi$, and its general form for arbitrary multiples of period vectors,  gives a surjection  $\pi :\Omega(\C) \to \Omega(\C')$.
(The details are given in the Appendix.) In particular while as a set  $\Omega(\C')$ can be "smaller" that
$\Omega(\C)$ (for example, horizontal lines with rational intercepts in $\Omega(\C)$ may be coalesced in  $ \Omega(\C')$ under $\pi$) the dimension of the spectrum (as indicated above) remains the same.

\medskip

\noindent {\bf Square-summable flexes.}
An infinitesimal flex being local represents the strongest form of rapid decay possible since it applies zero velocities to the framework points outside some bounded region.
It is natural to enquire to what extent a crystal framework $\C$ might be resistant to flexes whose velocities diminish to zero at infinity. With this
in mind  write $\K_a^2$ and $ \K_b^2$ for the Hilbert spaces of square summable sequences in
$\K_{atom}$ and  $\K_{bond}$.
Thus  $u=(u_{\kappa,k})\in \K^2_a$ is such that the sum of the squares of the Euclidean norms
$|u_{\kappa,k}|$ is finite.
It is elementary to show that $R(\C)$ then determines a bounded Hilbert space operator
from $\K_a^2$ to  $\K_b^2$. For a given translation group $\T$ this operator intertwines the associated shift transformations, as before, although now these transformations are unitary operators  on $\K_a^2$ and $ \K_b^2$. Identifying
square-summable sequences with square-integrable functions in a standard way one obtains unitary equivalences $U_a$ and $U_b$ between $\K_a^2$ and $L^2(\bT^d)\otimes \bC^{d|F_v|}$ and between  
 $\K_b^2$ and $L^2(\bT^d)\otimes \bC^{|F_e|}$ respectively. The corresponding
unitary transform $U_bR(C)U_a^*$ of the operator $R(\C)$  is then a multiplication operator between these matrix-valued function spaces and the
multiplying function is in fact the symbol function $\Phi_\C(\ol{z})$.
In this way the matrix function for $\C$ and its translation group appears naturally
from the point of view of square-summable velocity sequences.
For more details see Owen and Power
\cite{owe-pow-crystal} where other operator-theoretic considerations are
given.

More speculatively, it would be of interest to investigate other possible roles
of the matrix function, particularly with regard to approximate flexes and quantitative issues. For example for the 3D framework $\C$ we may define the non-negative scalar function
$\lambda$ on $\bT^3$ with
\[
\lambda : (z_1,z_2,z_3) \to \lambda_{\rm min}(\Phi_\C(z_1,z_2,z_3)^*\Phi_\C(z_1,z_2,z_3))
\]
where $\lambda_{\rm min}(A)$ denotes the smallest eigenvalue of the
positive operator $A$.
In particular, when the spectrum is trivial, that is, equal to the singleton set $\{(1,1,1)\}$  the function is nonvanishing except at this point and so $\lambda$
could be viewed as a measure of RUM resistance.

\section{RUMs and low energy phonons.}
In the less idealised setting of traditional mathematical crystallography,  mathematical models for crystalline dynamics assume that the atoms oscillate harmonically. The bond strengths are finite and a dynamical matrix
embodying them governs the modes and wave vectors of phonon excitations. We show how
the RUM spectrum $\Omega(\C)$ arises as the set of wave vectors ${\bf k}$ of
the harmonic excitations of $\C$ which induce vanishing bond distortion in their low frequency (long wavelength, low energy) limits.

Suppose that $\C$ is a crystal framework in $\bR^d$, with motif data $(F_v, F_e, \T)$
and suppose that the vertices $p_{\kappa, k}$ undergo
a  standard wave motion,
\[
p_{\kappa, k}(t) =  p_{\kappa, k} +u_{\kappa, k}(t),\quad \kappa\in F_v, k\in \bZ^d,
\]
where $u_{\kappa, k}(t)$ represents the local oscillatory motion of atom $\kappa$ in the translated unit cell with label $k\in \bZ^3$.
Following standard formula-simplifying conventions, the framework point motions  take values in $\bC^d$,  the case of real motion being recoverable from real and imaginary parts. (See Dove \cite{dov-book}.)
Thus it is assumed that we have
\[
u_{\kappa, k}(t)= u_\kappa\exp(2\pi i{\bf k}\cdot k)\exp(i\alpha t)
\]
where ${\bf u}=(u_\kappa)_{\kappa\in F_v}$ is a fixed vector in $\bC^{3|F_v|}$, where ${\bf k}$ is the wave vector and where $\alpha$ is the frequency.

Consider now the distortion $\Delta_e(t)$ for the edge $e=[p_{\kappa, k}, p_{\tau, k+\delta(e)}]$ measured as the change in the square of the edge length. We have
\[
\Delta_e(t):=|p_{\kappa, k}(t)- p_{\tau, k+\delta(e)}(t)|^2-|p_{\kappa, k}(0)- p_{\tau, k+\delta(e)}(0)|^2
\]
\[
=2Re\langle p_{\kappa, k}- p_{\tau, k+\delta(e)},u_{\kappa, k}(t)-u_{\tau, k+\delta(e)}(t)\rangle
\]
\[
+2Re\langle p_{\kappa, k}- p_{\tau, k+\delta(e)},u_{\kappa, k}(0)-u_{\tau, k+\delta(e)}(0)\rangle
\]
\[
 +
\epsilon({\bf u},{\bf k}, k,\alpha t)
\]
where
\[
\epsilon({\bf u},{\bf k}, k,\alpha t)=
|u_{\kappa, k}(t) - u_{k+\delta(e)}(t)|^2-|u_{\kappa, k}(0) - u_{k+\delta(e)}(0)|^2.
\]
First note that in any finite time period $[0,T]$ the difference quantities $\epsilon({\bf u},{\bf k}, k,\alpha t)$
tends to zero uniformly, for all $t\in [0,T]$ and all $k$ in $\bZ^3$, as the frequency
$\alpha$ tends to zero. This follows readily from the fact
that for any $\theta$ the quantity
\[
|\sin(\alpha t+\theta)-\sin (\alpha t) |^2 - |\sin (\theta)-\sin(0)|^2
\]
tends to zero uniformly for $t\in [0,T]$ as $\alpha$ tends to zero.

For the other terms
for $\Delta_e(t)$ note that
\[
2Re\langle p_{\kappa, k}- p_{\tau, k+\delta(e)},u_{\kappa, k}(t)-u_{\tau, k+\delta(e)}(t)\rangle
\]
\[
=2Re [e^{-\alpha t {-2\pi i{\bf k}}\cdot k}\langle p_{\kappa, k}- p_{\tau, k+\delta(e)},u_{\kappa}-\omega^{\delta(e)}u_{\tau}\rangle]
\]
\[
=2Re [e^{-\alpha t {-2\pi i{\bf k}}\cdot k}\langle p_{\kappa}- p_{\tau,\delta(e)},u_{\kappa}-\omega^{\delta(e)}u_{\tau}\rangle]
\]
which is zero, irrespective of $t$, if $(\omega^{k}u_{\kappa})$ is an infinitesimal flex of the framework.

It follows that we have proven the implication (i) implies (ii) in the following
proposition
and in fact the converse assertion follows from the same equations. The theorem underlies the correspondence of the points in  $\Omega(\C)$ with the wave vectors of RUM phonons
that arise in simulations.

\begin{theorem}\label{p:LongWavengthLimit}
Let $\C$ be a crystal framework, with specified periodicity,  and let ${\bf k}$ be a wave vector with point $\omega \in\bT^3$.
Then the following assertions are equivalent.

(i) $(\omega^ku_\kappa)_{\kappa,k}$ is a nonzero phase-periodic infinitesimal flex for $\C$.

(ii)
For the vertex wave motion
\[
p_{\kappa, k}(t) =  p_{\kappa, k} + u_\kappa\exp(2\pi i~{\bf k}\cdot k)\exp(i\alpha t),
\]
and a given time interval, $t\in [0,T]$,
the bond length changes
\[
\delta e(t)=|p_{\kappa, k}(t)- p_{\tau, k+\delta(e)}(t)|-|p_{\kappa, k}(0)- p_{\tau, k+\delta(e)}(0)|,
\]
tend to zero uniformly, in $t$ and $e$, as the wavelength $2\pi /\alpha$ tends to infinity.
\end{theorem}

In the last two decades the RUM spectra of  frameworks associated with specific material crystals have been derived by experiment and by simulation using lattice dynamics. Some of the results of this approach can be found in Giddy et al \cite{gid-et-al}, Hammond et al \cite{ham-dov-zeo1997}, \cite{hamdov-zeo1998}, Dove et al \cite{dov-exotic} and
Swainson and Dove \cite{swa-dov}.
In particular the programme CRUSH
has been used for this purpose and this method  reflects principle (ii) in the theorem above. Indeed in the simulations a double limiting process is used (the split atom method) in which each shared
vertex (often an oxygen atom) is duplicated, for each rigid unit, and connected by bonds of zero length and increasing strength, tending to infinity.
In this set up the RUM wave vectors coincide with those  for which the long wavelength limits have vanishing energy and through this connection  they can be identified in simulation experiments and counted.

\section{Determinations of RUM spectra.}
The rigid unit mode spectrum is now determined for a variety of basic crystal frameworks. Also we emphasise  an infinitesimal flex method for the identification of lines and planes of wave vectors. The spectrum is  of standard type (in the sense given in the remarks in Section 5) for
the frameworks $\C_{\bZ^d}, \C_{\rm sq}, \C_{\rm star},
\C_{\rm kag}, \C_{\rm Knet}, \C_{\rm Oct}$ and $ \C_{\rm SOD}$,
while for the 2D zeolite $\C_{\rm oct}$ it is a union of four closed curves.

Consider once again the basic grid framework $\C_{\bZ^2}$ in the plane with motif consisting of a single vertex, $F_v=\{p_\kappa\}$, and two edges. Examining all edges it becomes evident that there exists an  infinitesimal flex $u$  supported on the $x$-axis, with
 $u_{\kappa,(k_1,0)}=(1,0)$ for all $k_1\in\bZ$.
Using all the parallel translates $T_{(0,k_2)}u$ of $u$, we may define a phase-periodic
  velocity vector $v$ in $\K_{atom}$,
\[
v=\sum_{k_2\in\bZ}\omega_2^{k_2}T_{(0,k_2)}u,
\]
where $\omega_2$ is a fixed point in $\bT$. Note that $v$  is well-defined and
\[
R(\C)v = R(\C)\sum_{k_2\in\bZ}\omega_2^{k_2}T_{(0,k_2)}u=
\sum_{k_2\in\bZ}\omega_2^{k_2}R(\C)T_{(0,k_2)}u=\sum_{k_2\in\bZ}\omega_2^{k_2}T_{(0,k_2)}R(\C)u=0.
\]
Thus $v$ is an infinitesimal flex, phase-periodic for the point $(1,\omega_2)$ in $\bT^2$ and so $(1,\omega_2)$ lies in the RUM spectrum. In the language of wave vectors the RUM spectrum contains the line of wave vectors $(0,\alpha)$. By symmetry the line $(\alpha, 0)$ is also included.
Similar arguments apply to the kagome lattice which also has  linearly localised infinitesimal flexes. (See also \cite{gue-hut}, \cite{hut-fle} for example.)

More generally, suppose that a crystal framework $\C$  has a nonzero infinitesimal flex $u$ which is
\medskip

(a) \textit{band limited}, in the sense of being supported by a set of framework vertices within a finite distance of
a direction axis for $\T$, and

(b) \textit{periodic}, or more generally, \textit{phase-periodic} in the direction axis direction.
\medskip

By (a) one can form a sum analogous to that above, using the complementary axis direction(s), to obtain a well-defined phase-periodic infinitesimal displacement which,
by translational invariance and linearity, is an infinitesimal flex.
If  $\omega_1$ is the phase in (b)
then we deduce that $\{\omega_1\} \times \bT^{d-1}$ is contained in $\Omega(\C)$.

Thus, for the grid framework $\C_{\bZ^3}$ in three dimensions one deduces
from the evident line-localised infinitesimal flexes
that there are
three surfaces, $z=1$, $w=1$ and $u=1$, in $\Omega(\C_{\bZ^3})$.
In general a line-localised flex of this type leads directly to a hyperplane of wave vectors in the RUM spectrum.

Similar observations hold  for plane-localised
flexes. In three dimensions, for example, such a flex, which is assumed to be "in-plane phase-periodic", leads to a line of RUM wave vectors. This is the case for  $\C_{\rm Oct}$, considered below, and the RUM spectrum here is the union of these planes.
\medskip

\noindent {\bf Example (a): The regular 4-ring framework $\C_{\rm star}$.}
This 2D zeolite is defined by translates of the regular 4-ring of equilateral triangles in Figure  6.
It is sufficiently simple that
one can deduce its RUM spectrum and its crystal polynomial $p_{\rm star}(z,w)$  from infinitesimal arguments.

For a motif we may take $F_e$ to consist of the edges of the $12$-edged star and take $F_v$
to be the set of four vertices of the square together with the  westward and southern vertex. The four edges in the motif incident to the external vertices (north and eastward) provide four rows of the $12$ by $12$ matrix $\Phi_\C(z,w)$ each of which carries simple monomials (either $z$ or $w$ or their conjugates). Thus  $p_\C(z,w)$ has  total degree at most $4$. One can identify band-limited infinitesimal flexes  as indicated in Figures 6 and 7. Here the top and bottom vertices of each are fixed
and there is  horizontal periodic extension to a band-limited infinitesimal flex. In the former case there is two-step horizontal periodicity while in the latter case there is strict horizontal periodicity although the band is two cells wide.

\begin{center}
\begin{figure}[h]
\centering
\includegraphics[width=7.5cm]{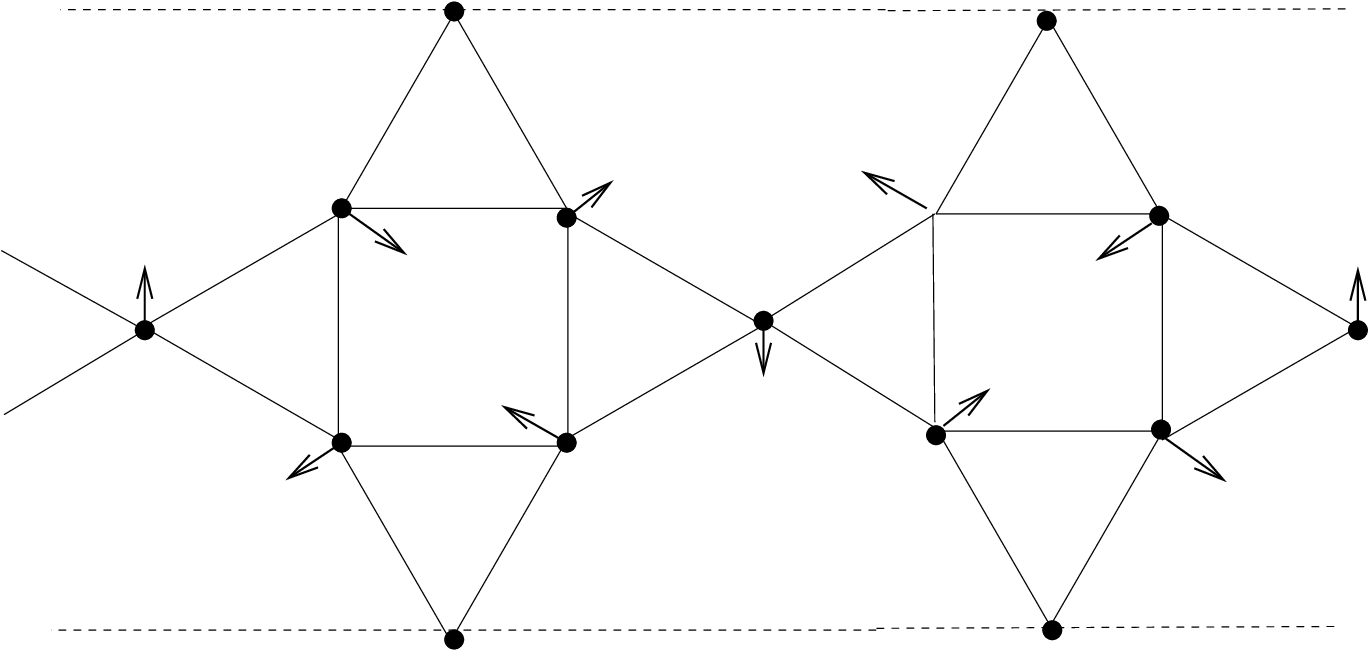}
\caption{A $2$-cell-periodic band-limited flex of $\C_{\rm star}$.}
\end{figure}
\end{center}

\begin{center}
\begin{figure}[h]
\centering
\includegraphics[width=5cm]{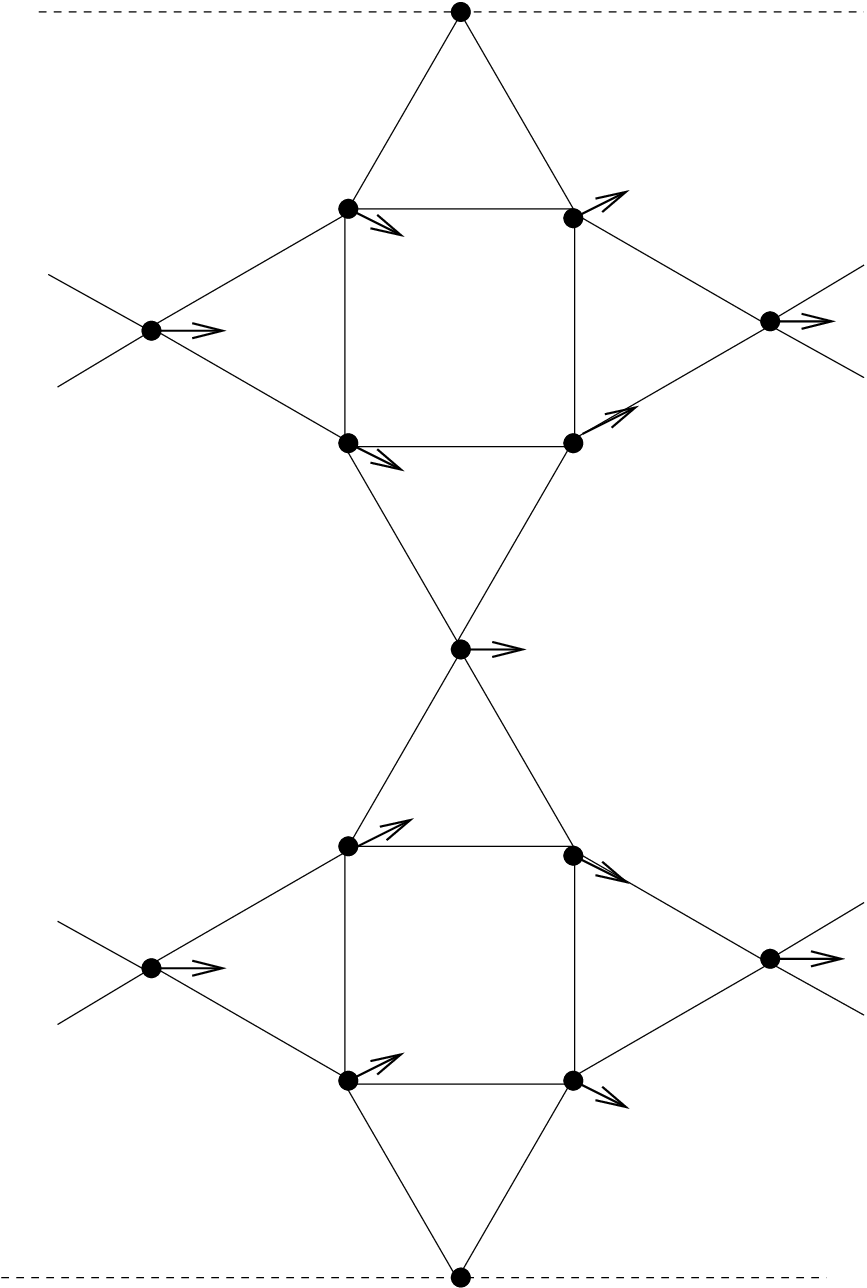}
\caption{A $1$-cell-periodic band-limited flex of $\C_{\rm star}$.}
\end{figure}
\end{center}

From the discussion above the first band-limited flex shows  that the phase $(-1,\omega_2)$ lies in  $\Omega(\C)$
for all $\omega_2 \in \bT$. By symmetry  $(\omega_1, -1)\in \Omega(\C)$
for all $\omega_1 \in \bT$.
The second band-limited flex shows that $\{1\}\times \bT$ lies in $\Omega(\C_{\rm star})$
and hence so too does $\bT\times \{1\}$ by symmetry.
Thus  $\Omega(\C)$ contains the set
\[
(\{1\}\times\bT) \cup (\bT \times\{1\}) \cup (\{-1\}\times\bT) \cup (\bT \times\{-1\})
\]
 and so $p_\C(z,w)$ must be divisible by the irreducible factors $z-1, w-1, z+1, w+1$. Since $p_\C(z,w)$ has total degree at most $4$ it follows that either $p$ vanishes identically or
\[
p_\C(z,w) = (z-1)(w-1)(z+1)(w+1).
\]
In fact the former case does not hold. One can see this, thematically, by demonstrating
that there are no local flexes or one may  compute  $\det \Phi(1/3,1/3) \neq 0$.
Thus the RUM spectrum is precisely the fourfold union above.

\medskip

\noindent {\bf Example (b): The 2D zeolite framework $\C_{\rm oct}$}.
There are no  local or  band-limited infinitesimal infinitesimal flexes  evident for the regular octagon framework and so the expectation  is that the RUM spectrum is trivial or a union of proper curves.

Returning to the 2D zeolites of Figure 5 the third of these, with external angle
$8\pi/12$, is  equal to $\C_{\rm star}$, although with a different translation group, $\T'$, for which the old period vectors are  rotated by $\pi/4$ and scaled by the factor $\sqrt{2}$.
In view of this rotation it follows  that
$$\Omega(\C_{\rm star},\T') = \{(w,w),(w,-w): w\in \bT\}.$$
In terms of reduced wave vectors this corresponds to the subset of the
unit square $[0,1)^2$ given as the union of the two diagonals.

As we have noted earlier, the 4-pointed star framework is related to its 8-pointed star companion $\C_{\rm oct}$ by a continuous flex. It follows that the $24$ by $24$ symbol matrix function $\Phi_{\rm star}(z_1,z_2)$ for the former (for $\T'$) is naturally "continuously connected" to the symbol function $\Phi_{\rm oct}(z_1,z_2)$ by an explicit continuous path $t\to \Phi_t(z_1,z_2)$. This in turn provides a set-valued map which we refer to as the \textit{RUM spectrum evolution} for this (periodicity-preserving) flex:
\[
t \to \Omega(\Phi_t(z)).
\]
When this is  made explicit by computation the octagon framework
has exotic (nonlinear) spectrum as indicated in Figure 8 and evolves towards the cross-shaped spectrum of $\C_{\rm star}$ under the continuous flex.

\begin{center}
\begin{figure}[h]
\centering
\includegraphics[width=4.5cm]{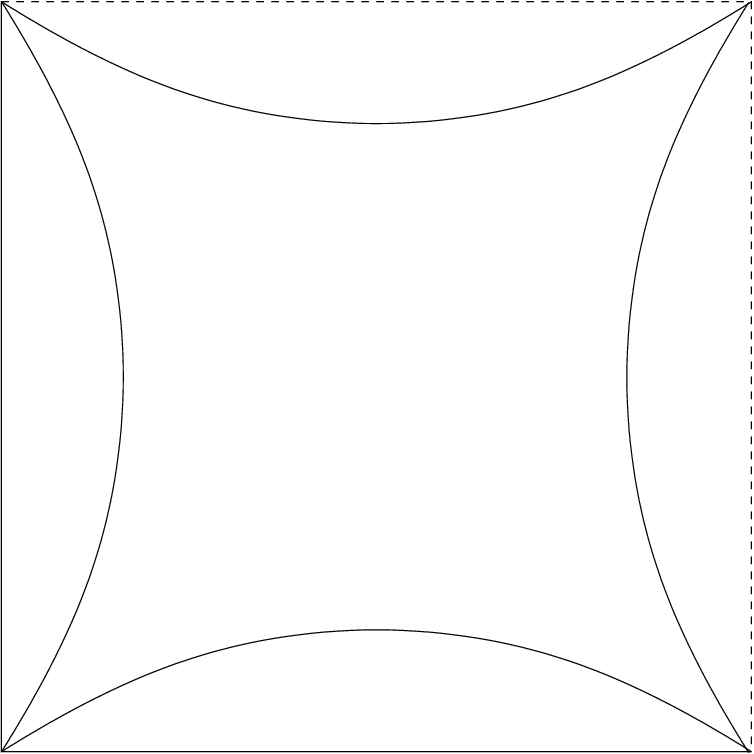}
\caption{The curved wave vector spectrum of $\C_{\rm oct}$.}
\end{figure}
\end{center}

In fact one can obtain the RUM spectrum of the octagon framework completely analytically,
although with some significant algebraic complexity, as follows.

Note first that a motif for $\C_{\rm oct}$
is formed by the 24 edges of the $8$-ring for $F_e$ with $F_v$ obtained by omitting four boundary vertices as, for example, in Figure 5. There are $8$ edges with external vertices
and each contributes a row to $\Phi_{\rm oct}$ with a simple monomial and so it follows that $p_{\rm oct}(z,w)$ has degree $8$ at most.
The $24 \times 24$ function matrix $\Phi_{\rm oct}(z)$  is sparse
and the (at most) four nonzero functions in each row may be conveniently normalised by dividing  by the magnitude of the $x$-coordinate
difference for that row. The magnitudes of the nonzero nonunit entries  are then the tangents of the angles $k\pi/24$, for $k=1,3,5,7,9,11$, all of which lie in the field extension $\bQ(\sqrt{2},\sqrt{3})$.
The crystal polynomial can be computed
and admits an explicit factorisation as the product
\[
p_{\rm oct}(z,w)=p_1(z,w)p_2(z,w),
\]
where
\[
p_1(z,w)=z^2w-(\sqrt{3}+\sqrt{2})zw^2+2(\sqrt{3}+\sqrt{2}-1)zw-(\sqrt{3}+\sqrt{2})z+w,
\]
\[
p_2(z,w)=z^2w-(\sqrt{3}-\sqrt{2})zw^2+2(\sqrt{3}-\sqrt{2}-1)zw-(\sqrt{3}-\sqrt{2})z+w.
\]
Each of the factors is responsible for two of the four closed curves that
comprise the RUM spectrum.

Returning to the as yet unconsidered 2D zeolite of Figure 5
(the first framework indicated, with an "8-ring of triangles encircling a square")
we remark that one can also show, by band-limited infinitesimal flex analysis,
that it  has  standard RUM spectrum, being the
subset of the unit square $[0,1)^2$ given as the union of the axes.

Each of these 2D zeolites has a 3D zeolite companion obtained by the layer construction.
The companion $\tilde{\C}_{\rm oct}$  for $\C_{\rm oct}$ also has exotic RUM spectrum and in fact by
earlier arguments contains the surface of points
$(z,w,u)$ in $\bT^3$ with $(z,u)$ in $\Omega(\C_{\rm oct})$ and $u$ any point of $\bT$.

As we have already remarked, the  two-dimensional crystal framework motion implied by Figure 5 is an example of a
\textit{finite flex} and  continuous and smooth flexes such as these serve to model  flexibility considerations for  zeolites and other micro-porous materials. These finite motions usually take place with an associated contraction and increase in rigid unit density. See for example
the collapsing mechanisms of Kapco et al \cite{kap-tre-tho-gue} and the
flexibility window determinations in Kapko et al \cite{kap-daw-tre-tho}.

\medskip

\noindent {\bf Example (c): A 2D zeolite with order N.}
Figure 9 shows a unit cell for a 2D zeolite, $\C_{\rm bowtie}$ say, which is of order $N$.
(This resolves an existence question posed by Mike Thorpe.) To see this property
one can verify that there is an infinitesimal flex
of the enclosed finite framework which assigns  zero velocities to the six boundary vertices and
a nonzero vertical velocity to the central vertex. Thus
the entire framework has a local infinitesimal flex and so the RUM spectrum is all of $\bT^2$.

We remark that in general it need not be the case that an order $N$ crystal framework
has a local infinitesimal flex internal to a unit cell.
For example, one could take a new motif and unit cell
in which the central vertex is shifted to the boundary and in this case one has to consider
a  threefold supercell before a local RUM appears.
\begin{center}
\begin{figure}[h]
\centering
\includegraphics[width=8cm]{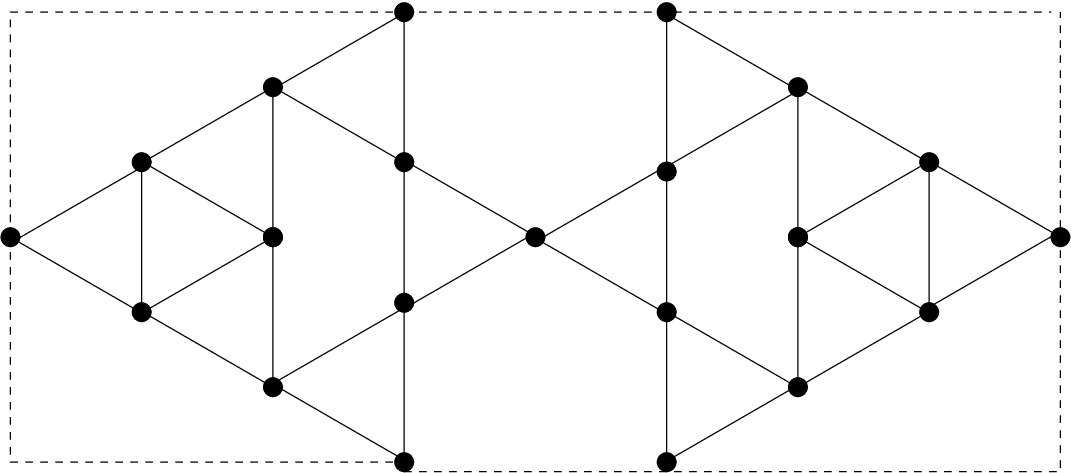}
\caption{A unit cell which defines $\C_{\rm bowtie}$.}
\end{figure}
\end{center}

\noindent {\bf Example (d): The kagome net framework  $\C_{\rm Knet}$ and its polynomial.}
The RUM spectrum of the 3D kagome net framework can be derived from that of the 2D kagome framework.
The spectrum of the latter is the zero set on the torus $\bT^2$ determined by
the crystal polynomial, which an earlier computation showed was equal to
$(z-1)(w-1)(z-w)$. One can derive  this from infinitesimal flex analysis as follows.
It is elementary to show that there is no local infinitesimal flex and so
$p_{\rm kag}(z,w)$ is necessarily nonzero. There are line-supported infinitesimal flexes in the directions of the period vectors $a_1$ and $a_2$
so it follows from the discussion above that $(z-1)$ and $(w-1)$ are factors. Let $u$ be a (similar) infinitesimal flex supported on a line in the direction $a_1-a_2$ and consider the infinitesimal flexes
\[
v=\sum_{k\in \bZ} \omega^kT_{(k,0)}u
\]
for $\omega \in \bT$. In view of the triangular symmetry of $\C_{\rm kag}$ in fact this flex is phase-periodic for the phase $(\omega,\omega)$ and it follows that $(z-w)$ is necessarily
a factor of $p_{\rm kag}(z,w)$. One can see, without calculation, that the total degree of this polynomial is at most $3$ and so the derivation is complete.

Moving up a dimension, a phase-periodic flex of $\C_{\rm kag}$, with phase $(\omega_1,\omega_2)$ say, induces a "layer-limited" infinitesimal  flex of the kagome net framework $\C_{\rm Knet}$. Thus, for all $\omega \in \bT$ there is an infinitesimal flex of $\C_{\rm Knet}$ with phase $(\omega_1,\omega_2,\omega)$. Similar assertions hold for the other two translation group planes.
The crystal polynomial $p_{\rm Knet}(z,w,u)$ has total degree at most $6$ and  must vanish on the six planes
$z-1=0, w-1 =0,u-1=0, z-w=0, w-u=0, z-u=0$. It follows that
\[
p_{Knet}(z,w,u) = (z-1)(w-1)(u-1)(z-w)(w-u)(z-u),
\]
for the monomial order with  $z>w>u$.

The polynomial above was also obtained in Wegner \cite{weg} and Owen and Power \cite{owe-pow-crystal} by direct calculation.

The kagome net framework and its various periodic positions or placements feature as the
tetrahedral net frameworks for a range of materials and their phases. It is the framework for $\beta$-cristobalite, for example,  while a particular placement gives the framework for tridymite. This
was the first material for which curved surfaces of RUMs were observed.
(Dove et al \cite{dov-exotic}.)

\medskip

\noindent {\bf Example (e): Sodalite and $\C_{\rm SOD}$.}
The framework $\C_{\rm SOD}$ has a symbol function with $72$ rows and columns. Indeed, it is in Maxwell counting equilibrium, being a 3D zeolite crystal framework, and the motif set $F_e$ consists of the edges of three 4-rings of tetrahedra.
We prove  that $\C_{{\rm SOD}}$ is of order $N$. Specifically we show, by infinitesimal flex geometry, that there is a nonzero  infinitesimal flex $v$ of the finite sodalite cage framework such that all the outer vertices are fixed by $v$.  That is,  $v_{\kappa, \delta}=0$ if $p_{\kappa,\delta}$ is any of the 24  outer vertices of the cage. The "outer fixed" sodalite cage framework has $36$ free vertices with $108$ degrees of freedom while there are $144$ constraining edges. Despite this considerable over-constraint there is sufficient symmetry to allow in a proper infinitesimal flex.

We shall show  that an individual  $4$-ring, $R_1$ say, of the sodalite cage
has an infinitesimal flex, $v^{(1)}$ say, which fixes a coplanar quadruple of "outer" vertices, such as the upper vertices of the $4$-ring in Figure 4,
and flexes the other quadruple in their common plane.
These four velocites have equal magnitude and
in Figure 4 are directed towards two opposing corners of the imaginary cube. (See the flex arrows in Figure 4). Taking  $v^{(1)}$ so that these vectors have magnitude $1$ it follows that
$v^{(1)}$ is determined up to sign and that this sign may be specified by labeling the cube corners "a" and "r" for their attracting and repelling sense. Note that one can label the eight corners of the imaginary cube in this manner so that no like labels are adjacent. In this case the individual flexes  $v^{(1)}, \dots , v^{(6)}$ of the six $4$ rings of the sodalite cage have equal displacement vectors at common vertices.
This consistency shows that there is an infinitesimal flex of the entire sodalite cage
in which the outer vertices are fixed, as required.

It remains to show that there is the stated flex of the $4$-ring $R_1$. To this end
let $p_1, p_2$ be two non-opposite top vertices of $R_1$ with intermediate vertex $p_3$, let $p_1, p_3, p_4$ be the vertices of an inward facing face of a tetrahedron of $R_1$ with vertices
$p_1, p_3, p_4, p_5$, so that the lower vertex $p_5$ is a cube-edge midpoint.
There is a unique "inward and upward" displacement velocity $u_3$ of the intermediate vertex $p_3$ which has unit length and is such $(u_1,u_2,u_3)=(0,0,u_3)$ is a flex for the two edges $[p_1,p_3], [p_3,p_2]$. The displacement vector $u_3$ induces a unique displacement vector $u_5$ which is in the direction of the cube edge and is such that
\[
\langle u_5-u_3,p_5-p_3 \rangle =0.
\]
The triple $u_1=0, u_3$ and $u_5$ now determine the infinitesimal motion of the tetrahedron,
with flex vector $u_4$ for $p_4$. However, the reversal (sign change) of $u_3$ induces the
reversal of $u_5$ so it is clear from the symmetric position of the tetrahedron that
$u_4$ must be the unique unit norm "outward and downward" flex at $p_4$. Continuing around the ring it follows that $R_1$ has the desired infinitesimal flex.

One can apply similar constructive flex arguments to other zeolite frameworks and of course to any zeolite crystal framework
which contains a sodalite cage as above, such as $\C_{\rm LTA}$.
Also we note (as do Kapko et al \cite{kap-daw-tre-tho}) that  $\C_{{\rm RWY}}$ is derived from $\C_{{\rm SOD}}$ by replacing each
tetrahedron by a rigid unit of four tetrahedra. Thus the same infinitesimal flex geometry
applies and  $\C_{{\rm RWY}}$ has order $N$.

\medskip

\noindent {\bf Example (f): Perovskite, $\C_{\rm sq}$ and  $\C_{\rm Oct}$.}
Consider the integer translation group
$\T$ and the determination of the framework through the primitive motif $(F_v, F_e)$ where
\[
F_v=\{0,1/2,1/2), (1/2,0,1/2),(1/2,1/2,0)\} =\{p_{\kappa,0}: 1\leq \kappa \leq 3\}
\]
and where $F_e$ consists of the twelve framework edges between the centres of adjacent faces
of the unit cube $[0,1]^3$. Thus the vertices of $V(F_e)\backslash F_v$ have the form $p_{i,\gamma_i}$, where
$\gamma_1= (1,0,0), \gamma_2=(0,1,0), \gamma_3=(0,0,1)$.
The framework is therefore edge rich and the matrix function $\Phi_{\rm Oct}(z)$ is 12 by 9.

The framework $\C_{\rm Oct}$ is a 3D analogue of the 2D squares framework $\C_{\rm sq}$
and may be obtained from it by a layer construction and a discarding of redundant edges internal to the octahedra. Thus infinitesimal flexes of the 2D squares lattice imply plane-localised flexes for $\C_{\rm Oct}$. This observation can be made the basis for an infinitesimal flex analysis determination of the RUM spectrum. A more algebraic approach
is possible as follows.

Performing row operations on $\Phi_{\C_{\rm sq}}(z)$, as given in Section 2, we see that
$(\ol{z},\ol{w})$ is a point of the RUM spectrum if and only if the equivalent matrix
{\small
$$ \Psi(z,w)= \left[ \begin {array}{cccc}
1 &-1&-1&1\\
0 &-2&z-1&z+1\\
0 &0&1-wz&1-wz\\
0 &0&0&-2z+2w\\
0&0&-2+2z&0
\end {array} \right]
$$
}
has rank less than $4$. This occurs if and only if
{\small
$$\left[ \begin {array}{cc}
1-wz&1-wz\\
0&-2z+2w\\
-2+2z&0
\end {array} \right]
$$
}
has rank equal to $0$ or $1$. The rank is $0$ if and only if $z=w=1$,
corresponding to the two-dimensional space of rigid motions
with phase $(1,1)$, and the rank is $1$ if and only if $z=w=-1$.
Thus
\[
\Omega(\C_{\rm sq})=\{(1,1),(-1,-1)\}.
\]
The infinitesimal flex for the phase $(-1,-1)$ is the one for which the rigid units,
in this case squares with diagonals, rotate infinitesimally in alternating senses.

The alternating rotation flex of $\C_{\rm sq}$ induces a
plane-localised flex of $\C_{\rm Oct}$
in each of the framework planes $x=1/2, y=1/2, z=1/2$.
It follows that $\Omega(\C_{\rm Oct})$  contains the three sets of phases,
\[
\bT\times\{-1\}\times\{-1\}, \quad \{-1\}\times \bT\times\{-1\}. \quad \{-1\}\times \{-1\}\times \bT.
\]
That the spectrum is no more than the union of these sets and the singleton $(1,1,1)$ can be seen from a row analysis of the $12$ by $9$ function matrix $\Phi_{\C_{\rm Oct}}(z)$
in the same style as the argument for $\C_{\rm sq}$.
Thus, in wave vector formalism, the RUM spectrum of the octahedral net
 $\C_{\rm Oct}$ is the set of lines
\[
(\alpha,{1}/{2},{1}/{2}),\quad ({1}/{2},\alpha,{1}/{2}),\quad ({1}/{2},{1}/{2},\alpha)
\]
together with the  wave vector $(0,0,0)$.

The corner connected octahedron net crystal  framework $\C_{\rm Oct}$ is associated with cubic perovskites, such as ~SiTO$_3$, and RUM distributions have been determined experimentally, Giddy et al \cite{gid-et-al}, Dove et al \cite{dov-exotic}.

\section{Appendix}

\noindent {\bf The periodic floppy mode counting formula.}
For notational clarity we assume that $d=3$.
Let $r=(r_1,r_2,r_3)$ and consider the finite-dimensional space  $\K_{a}^r$
of $r$-periodic complex velocity vectors.

Write $k\in r$ to denote $k=(k_1,k_2,k_3) $ with $0\leq k_i < r_i$.
For $\omega_l=e^{2\pi i/{r_l}}$, for $l=1,2,3$,  write $\bT^3_r $ for  "discrete torus"
\[
\bT^3_r=\{\omega=(\omega_1^{k_1}, \omega_2^{k_2}, \omega_3^{k_3}):k\in r\}.
\]
If $z=(z_1,z_2,z_3)$ is a point of the usual $3$-torus $\bT^3$  we  write
$z^k$ for the product $z_1^{k_1}z_2^{k_2}z_3^{k_3}$ in $\bT^3$. Similarly,
with $W_1, W_2, W_3$ defined as the shift transformations 
$T_{\gamma_1}, T_{\gamma_2}, T_{\gamma_3}$ restricted to the space 
$\K_a^r$
we write $W^k$ for the product
$W_1^{k_1}W_2^{k_2}W_3^{k_3}$. In particular $W^r=I$.

Note that if $u$ is an $r$-fold periodic flex then the velocity vector
\[
u'=\sum_{k\in r} W^ku
\]
is \emph{strictly} periodic. Since $R(\C)$ commutes with the shifts the velocity vector  $u'$ is a sum of infinitesimal
flexes and so is a strictly periodic infinitesimal flex.

Similarly, if $\omega \in \bT^3_r $ then
\[
u_\omega=\sum_{k\in r} \omega^kW^ku
\]
is an infinitesimal flex which is phase-periodic for $\ol{\omega}$. Since we have the recovery formula
\[
u = \frac{1}{r_1r_2r_3}\sum_{\omega \in \bT^3_r} u_\omega
\]
it follows that the space of $r$-fold periodic infinitesimal flexes is the
direct sum of the space of $\omega$-phase periodic infinitesimal flexes.
The counting formula now follows.
\bigskip

\noindent{\bf Surjectivity of $\pi: \Omega(\C) \to \Omega(\C')$.}
Similarly, let $\T'$ be the subgroup of the translation group $\T=\{T_k:k\in \bZ^3\}$ for $\C$ which is
associated with $(r_1,r_2,r_3)$, let $\C'$ be the associated crystal framework and suppose that $u$ is a nonzero  phase-periodic infinitesimal flex of $\C'$ with multi-phase $\eta$ in $\bT^3$.

Let $\bT^3_{r,\eta}$ be the set of points $(z_1,z_2,z_3)$ where $z_i$ ranges over the $r_i$ roots of $\eta_i$.
If $\omega \in \bT^3_{r,\eta} $ then the velocity vector
\[
u_\omega=\sum_{k\in r} \omega^kW^ku
\]
is an infinitesimal flex for $\C$ which is  phase-periodic for $\C$, with multi-phase $\ol{\omega}$.
Also we have the recovery formula
\[
(\eta_1\eta_2\eta_3)u = \frac{1}{r_1r_2r_3}\sum_{\omega \in \bT^3_r} u_\omega.
\]
It follows that at least one of the flexes $u_\omega$ is nonzero. That is, there is an $r$-fold root of $\eta$ in the RUM spectrum
$\Omega(\C)$ and the surjectivity of the map $\pi$ follows.

\bigskip
\noindent{\bf The existence of local infinitesimal flexes.} 
We prove the following.

\begin{theorem}\label{p:flopexists}Let $\C$ be a crystal framework in $\bR^d$ which is of maximal order ("order $N$") for periodic floppy modes. Then $\C$ has a local infinitesimal flex.
\end{theorem}

\begin{proof}
Consider the vector space $\V_n$ say of $n$-fold periodic real velocity vectors. This has dimension $|dF_v|n^d$ as $n$ goes to infinity. By assumption the subspaces $\ker R_n(\C)$ have dimensions of order $n^d$. Fix the motif $M=(F_v,F_e)$ for $\C$ and consider the natural motifs $M_n$ for the $n$-fold translation group which are formed by translates of $M$ ($n^d$ translates in fact).
(One could arrange $M_n \subseteq M_{n+1}$ but this is not necessary for the argument.) The motif $M_n$, which is a pair $(F_v(n),F_e(n))$, has \emph{boundary vertices} by which we mean the vertices of edges in $F_e(n)$ which are not vertices in $F_v(n)$. Note that the cardinality of these sets gives a sequence of order $n^{d-1}$. Thus the vector subspace, $\B_n$ say,  
of $n$-fold periodic real velocity vectors which assign zero velocities
to the nonboundary vertices has dimension of order $n^{d-1}$.

Let 
\[
P_n:\V_n \to \B_n
\]
be linear transformations that are projections. Then, in view of the
order of dimension growth elementary linear algebra shows that there is a nonzero vector $u$ in $\ker R_n(\C)$ for some large enough $n$ such that
$P_n(u)=0$. 

Let $u'$ be the velocity vector which agrees with $u$ for components for the non boundary framework points of $M_n$ and is defined to be zero for all other coordinates. 

Since the $n$-fold periodic flex $u$ "vanishes on the boundary of the $n$-fold supercell" in the sense above one can readily check that $u'$ is an infinitesimal flex of $\C$. Also $u'$ is nonzero and finitely supported, as desired.
\end{proof}


{\bf Acknowledement.} The development
relates to themes and problems outlined at a  London Mathematical Society workshop
on \textit{The Rigidity of Frameworks and Applications} held at Lancaster University in July 2010. The author is grateful for discussions and communications with Martin Dove, Simon Guest, John Owen, Mike Thorpe, Franz Wegner and Walter Whiteley. The project is supported by
the EPSRC grant  EP/J008648/1 on "Crystal Frameworks, Operator Theory and Combinatorics".



\bibliographystyle{abbrv}
\def\lfhook#1{\setbox0=\hbox{#1}{\ooalign{\hidewidth
  \lower1.5ex\hbox{'}\hidewidth\crcr\unhbox0}}}

\end{document}